\documentclass[a4paper,12pt,reqno]{amsart}
\usepackage{paralist} 
\usepackage[T1]{fontenc}
\usepackage{geometry}
\allowdisplaybreaks

\usepackage{amsmath}
\allowdisplaybreaks
\usepackage{amssymb}
\usepackage{amsfonts}
\usepackage{color}
\usepackage{enumitem}

\usepackage{booktabs}

\usepackage{etoolbox}
\patchcmd{\thebibliography}{\section*{\refname}}{}{}{}

\usepackage{verbatim}

\usepackage{graphicx}
\usepackage{subfig}
\usepackage{bm}
\usepackage{tikz}
\usepackage{cite}
\bgroup

\numberwithin{equation}{section}

\usepackage[bookmarks,bookmarksopen=false,% Klappt die Bookmarks in Acrobat aus
                  pdfauthor={Autor},%
                  pdftitle={Titel},%
                  colorlinks=true,%
                  linkcolor=blue,%
                  citecolor=blue,%
                  urlcolor=blue,%
                ]{hyperref}

\newtheorem{lemma}{Lemma}[section]
\newtheorem{theorem}[lemma]{Theorem}
\newtheorem{corollary}[lemma]{Corollary}
\newtheorem{proposition}[lemma]{Proposition}

\newtheorem{conj.}[lemma]{Conjecture}
\newtheorem{definition}[lemma]{Definition}

\newtheorem{remark}[lemma]{Remark}

\newcommand{\pbox}{\hfill$\Box$\\}

\newcommand{\R}{\mathbb{R}}

\renewcommand{\le}{\leqslant}
\newcommand{\C}{\mathbb{C}}

\newcommand{\HS}{\mathcal{HS}}

\renewcommand{\H}{\mathcal{H}}

\newcommand{\lf}{localised }

\definecolor{darkviolet}{rgb}{0.58,0,0.83} %{148,0,211}

\usetikzlibrary{patterns}

\title[Kernel theorems for   localised frames]{Kernel theorems for 
  operators on
co-orbit spaces 
associated with localised frames}

\author{Dimitri Bytchenkoff}
\address{(D. B.) Acoustics Research Institute, Austrian Academy of Sciences, Wohlle- bengasse 12-14, 1040 Vienna, Austria;
\newline
\indent Faculty of Mathematics, University of Vienna, 
Oskar-Morgenstern-Platz 1, 1090 Vienna, Austria;
\newline 
\indent Université de Lorraine, Laboratoire d’Energétique et de Mécanique Théorique et Appliquée, 2 avenue de la Forêt de
Haye, 54505 Vandoeuvre-lès-Nancy, France}
\email{dimitri.bytchenkoff@oeaw.ac.at; dimitri.bytchenkoff@univie.ac.at; dimitri.bytchenkoff@univ-lorraine.fr}
\author{Michael Speckbacher}
\address{(M. S.) Faculty of Mathematics, 
University of Vienna,
Oskar-Morgenstern-Platz 1, 
A-1090 Vienna, Austria}
\email{michael.speckbacher@univie.ac.at}
\author{Peter Balazs}
\address{(P. B.) Acoustics Research Institute, Austrian Academy of Sciences, Wohlle- bengasse 12-14, 1040 Vienna, Austria }
\email{peter.balazs@oeaw.ac.at}
\date{}

\begin{document}

\begin{abstract}
 Kernel theorems  provide a convenient representation of bounded linear operators. For the operator acting on a concrete function space, this means that its action on any element of the space can be expressed as a generalised integral operator, in a way reminiscent of the matrix representation of linear operators acting on finite dimensional vector spaces. We prove kernel theorems for bounded linear operators acting on co-orbit spaces associated with \lf frames. Our two main results characterise the spaces of operators whose generalised integral kernels belong to the co-orbit spaces of test functions and distributions  associated with the tensor product of the \lf frames respectively.  Moreover, using a version of Schur's test, we establish a characterisation of the bounded linear operators between some specific co-orbit spaces and kernels in mixed-norm co-orbit spaces.
\end{abstract}

\maketitle

\noindent\textit{2020 Mathematics subject classification:} {42B35, 42C15, 46A32, 47B34}\\
\noindent \textit{Keywords:} {Kernel theorem, localised frames, co-orbit spaces, operator representation, tensor products}

\section{Introduction}

\noindent 
This paper is  concerned with co-orbit spaces, which were originally introduced in the late 1980's by  Feichtinger and  Gröchenig %in their celebrated co-orbit theory 
\cite{Feichtinger_1988, Feichtinger_1989,Feichtinger_1989_2}, and the bounded linear operators that act on these spaces. Co-orbit spaces, such as modulation \cite{Feichtinger_1981} and Besov spaces \cite{Besov_1961}, are defined as   Banach spaces of functions or distributions whose membership is determined by integrability conditions imposed on their generalised wavelet transform. Co-orbit theory provides a convenient way of analysing and synthesising spaces of functions or distributions. The original definition of co-orbit spaces relied on an integrable representation of a locally compact group on a Hilbert space. In the early 2000's,  Gröchenig and  Fornasier showed how \lf frames can be used to introduce a novel class of co-orbit spaces and guarantee their stable decomposition  \cite{Fornasier_2005, Groechenig_2003, Groechenig_2004}. In this report, we prove kernel theorems for bounded linear operators acting on co-orbit spaces associated with \lf frames.

The main purpose of any kernel theorem is to {express} a given operator as a generalised integral operator. The original kernel theorem, proved by Laurent Schwartz in 1952, states that for any continuous linear operator $O$ mapping the Schwartz space $S (\mathbb{R}^n)$ on the space of tempered distributions $S' (\mathbb{R}^n)$ there is   a unique tempered distribution $K \in S' (\mathbb{R}^{2n})$, called the \emph{kernel}, that satisfies
\begin{equation}
\langle O f_1 ,  f_2 \rangle_{S' (\mathbb{R}^n),  S (\mathbb{R}^n)}
= \langle K ,   f_1 \otimes f_2 \rangle_{S' (\mathbb{R}^{2n}),   S (\mathbb{R}^{2n}) },
\label{eq:Schwartz}
\end{equation} 
for any $f_1,f_2 \in S (\mathbb{R}^n)$. Furthermore, if $K$ happens to be a locally integrable function, then
\begin{equation}
\langle O f_1,  f_2 \rangle_{S' (\mathbb{R}^n),  S (\mathbb{R}^n)}
= \int_{\mathbb{R}^{n}} \int_{\mathbb{R}^{n}} K (x, y) f_1 (y) \overline{f_2 (x) }dy dx,
\label{eq:Schwartz_2}
\end{equation} and so the operator $O$ is indeed expressed as an integral %or, so to speak, an integral
operator. Other classical kernel theorems were established, for example, for the continuous operators mapping $\mathcal{D}(\R^d)$ on $\mathcal{D}'(\R^d)$ \cite[Theorem~5.2]{testfunction}, and for operators mapping Gelfand-Shilov spaces on their distribution spaces \cite{geshi}. 

In 1980, Feichtinger proved a kernel theorem for bounded linear operators mapping the modulation spaces $M^1(\mathbb{R}^n)$ to $M^\infty(\mathbb{R}^n)$ \cite{Feichtinger_1980}, the spaces of choice in time-frequency analysis \cite{Feichtinger_1981}. Using a version of Schur's test, Cordero and Nicola characterised the kernel of operators mapping $M^1(\mathbb{R}^n)$ on $M^p(\mathbb{R}^n)$ as well as those mapping $M^p(\mathbb{R}^n)$ on $M^\infty(\mathbb{R}^n)$ \cite{Cordero_2019}. Finally, Balazs, Gr\"ochenig and Speckbacher took a more abstract approach to prove their kernel theorem \cite{Balazs_2019} for operators acting on co-orbit spaces associated with integrable representations of locally compact groups. 
Indeed, they considered irreducible unitary integrable representations $\pi_n:G_n\to \H_n$, $n=1,2,$   of locally compact groups $G_n$   that generate two scales of co-orbit spaces and showed a one-to-one correspondence between   bounded linear operators mapping $\text{Co} L^1(G_1)$ to $\text{Co} L^\infty(G_2)$ and their   kernel   belonging to the co-orbit space $\text{Co} L^\infty(G_1\times G_2)$ generated by the tensor product representation $\pi_1 \otimes \pi_2:G_1 \times G_2\to \H_1\otimes \H_2$. 

In this paper, we go even further and prove kernel theorems for bounded linear operators acting on co-orbit spaces generated by \lf frames  \cite{Fornasier_2005,Fornasier_Rauhut} and, by doing so, show that the group structure and the corresponding representation theory are an unnecessary, albeit very elegant setup for proving the kernel theorem. {Here is the essence of our main results. See Sections~\ref{sec:background}
and \ref{sec:tensor-of-frames}
for the definitions of the involved co-orbit spaces and Section~\ref{sec:proofs} for detailed statements and proofs of the subsequent results. }
\medskip

\noindent\textbf{Theorem.} (Outer kernel theorem) \textit{ Let $\Psi_n$ be $\mathcal{A}_n$-localised frames and $w_n$ an admissible weight for $n=1,2$. A linear operator $O$ from $\H^1_{w_1}(\Psi_1)$ on $\H_{1/w_2}^\infty(\Psi_2)$ is bounded if and only if there exists a unique kernel $K\in \H^{\infty,\infty}_{1/(w_1\otimes w_2)}(\Psi_1\otimes \Psi_2)$ such that 
\begin{equation}
\langle Of_1,f_2\rangle_{\H^\infty_{1/w_2}(\Psi_2),\H^1_{w_2}(\Psi_2)}=\langle K,f_1\otimes f_2\rangle_{\H^{\infty,\infty}_{1/(w_1\otimes w_2)}(\Psi_1\otimes\Psi_2),\H^{1,1}_{w_1\otimes w_2}(\Psi_1\otimes \Psi_2)},
\end{equation}
for all $f_1\in \H^1_{w_1}(\Psi_1)$ and $f_2\in \H^1_{w_2}(\Psi_2).$}
\medskip
 
\begin{comment}
Finally, we prove   our version of the \emph{inner kernel theorem} \cite{Feichtinger_Jakobsen} describing the space of operators with kernel in $\H^{1,1}_{w_1\otimes w_2}(\Psi_1\otimes\Psi_2)$, see Theorem~\ref{innere_kernsatz}.
\end{comment}

\noindent\textbf{Theorem.} (Inner Kernel Theorem)
\textit{Let $\Psi_n\subset \H_n$ be an $\mathcal{A}_n$-localised frame, and $w_n$ an $\mathcal{A}_n$-admissible weight for $n=1,2$. 
 There is an isomorphism between kernels in $\H^{1,1}_{w_1\otimes w_2}(\Psi_1\otimes \Psi_2)$ and the space
\begin{align*}
\mathcal{B}:=  \Big\lbrace &O\in  B\big(\H^\infty_{1/w_1}(\Psi_1),\H^1_{w_2}(\Psi_2)\big):\ \text{there exist } f_{  r} \in H^1_{w_1}(\Psi_1),g_r\in H^1_{w_2}(\Psi_2)  \text{ s.t.}\\ & O  =  \sum_{r \in R}\langle\, \cdot\,  ,f_{ r}\rangle_{H^\infty_{1/w_1}(\Psi_1),H^1_{w_1}(\Psi_1)} \,  g_{  r}   \text{ and }   \, \sum_{r \in R}  \| f_{ r} \|_{H^1_{w_1}(\Psi_1)} \,  \| g_{   r} \|_{H^1_{w_1}(\Psi_1)} < \infty  \Big\rbrace.
\end{align*}}

The terminology of  outer and inner kernel    theorems  \cite{Feichtinger_Jakobsen} originates in the fact that the corresponding kernels belong, respectively, to the outer- and innermost spaces on the scale of co-orbit spaces, as can be seen in the following chain of topological inclusions
$$ \H^{1,1}_{w_1\otimes w_2}(\Psi_1\otimes\Psi_2)\hspace{-.5pt}\subset\hspace{-.5pt} \H^{p,p}_{w_1\otimes w_2}(\Psi_1\otimes\Psi_2)\hspace{-.5pt}\subset\hspace{-.5pt} \H^{p,p}_{1/(w_1\otimes w_2)}(\Psi_1\otimes\Psi_2)\hspace{-.5pt}\subset \hspace{-.5pt}\H^{\infty,\infty}_{1/(w_1\otimes w_2)}(\Psi_1\otimes\Psi_2),$$ for weights $w_1,w_2 \geq 1$, and $p \in (1, \, \infty).$

Following the ideas developed in \cite{Cordero_2019}, we apply a version of Schur's test to characterise the kernels of the bounded operators mapping $\H^1_{w_1}(\Psi_1)$ to $\H^p_{w_2}(\Psi_2)$ and those mapping $\H^p_{w_1}(\Psi_1)$ to $\H^\infty_{w_2}(\Psi_2)$, see Theorem~\ref{Supergeiler_Satz}.

\begin{comment}

\textcolor{purple}{For weights $w_n\geq 1, n=1,2,$ and $1<p<\infty,$ the inclusions $$ \H^{1,1}_{w_1\otimes w_2}(\Psi_1\otimes\Psi_2)\hspace{-.5pt}\subset\hspace{-.5pt} \H^{p,p}_{w_1\otimes w_2}(\Psi_1\otimes\Psi_2)\hspace{-.5pt}\subset\hspace{-.5pt} \H^{p,p}_{1/(w_1\otimes w_2)}(\Psi_1\otimes\Psi_2)\hspace{-.5pt}\subset \hspace{-.5pt}\H^{\infty,\infty}_{1/(w_1\otimes w_2)}(\Psi_1\otimes\Psi_2),$$
motivate the notion of inner and  outer kernel theorem 
as the corresponding spaces of kernels are the inner- and outermost spaces in this chain of Banach spaces.
}

We also characterise the bounded operators from $\H^1_{w_1}(\Psi_1)$ to $\H^p_{w_2}(\Psi_2)$ and those from $\H^p_{w_1}(\Psi_1)$ to $\H^\infty_{w_2}(\Psi_2)$ following the ideas of \cite{Cordero_2019}, see Theorem~\ref{Supergeiler_Satz}. 
\end{comment}

%A few words about what motivated us to undertake this whole investigation. 
Let us shortly motivate to study kernel theorems within this setup.
There is every reason to use specifically localised frames here.
Firstly, \lf frames allow for essentially sparse atomic decompositions of the elements of their associated co-orbit spaces. Secondly, the Gram matrices of \lf frames, by  definition, are essentially sparsely populated. Therefore, using \lf frames, together with Galerkin's method \cite{Galerkin_1915, Balazs_2008}, can be expected to produce sparse matrices that represent the bounded linear operators that act on  co-orbit spaces    \cite{Futamura_2009}. This in itself can, as highlighted in \cite{Balazs_2008, Balazs_2017}, be very useful for solving corresponding operator equations numerically. Moreover, as we shall show in this paper, the kernel of any bounded linear operator mapping one co-orbit space on the other can be synthesised from the matrix that represents the operator, therefore leading to a sparse decomposition of the kernel.

This paper is organised as follows. In Section~\ref{sec:background} we provide the necessary frame theoretical background, namely the basics of the theory of \lf frames. In  Section~\ref{sec:tensor-of-frames} we {show that co-orbit spaces generated by tensor products of \lf frames share some properties with those generated by \lf frames.} In Section~\ref{sec:proofs} we state and prove our main results, namely we prove the outer and inner kernel theorems for co-orbit spaces. % in its most general form,state and prove the inner kernel theorem, make
Moreover, we prove statements about projective tensor products of two co-orbit spaces, boundedness of linear operators between certain intermediate co-orbit spaces %co-orbit spaces associated with a tensor product frame 
 and provide easily verifiable conditions under which a linear operator %mapping Hilbert spaces on one another
 belongs to the Schatten-$p$ classes.

Note that in \cite{sampta} the same authors stated the outer kernel theorem in a simplified form and gave a sketch of its proof.

\section{Background information and notation}\label{sec:background}

\noindent  %{In this section we recall the notions and facts that we will use to prove our kernel theorem.}
{Throughout this paper we will write $A\lesssim B$ if there exists $C\geq 1$ such that $A\leq CB$ and  $A\asymp B$ if $C^{-1}B\leq A\leq CB$. Moreover, $\H,\H_1,\H_2$ will always denote separable Hilbert spaces.} We shall use brackets $\langle \hspace{1pt}  \cdot\hspace{1pt} ,\hspace{-0.5pt}\cdot\hspace{1pt}  \rangle_{\H}$ for the scalar product in the Hilbert space and  $\langle \hspace{1pt}  \cdot\hspace{1pt} ,\hspace{-0.5pt} \cdot\hspace{1pt}  \rangle_{X,X^\ast}$ for the scalar product of the element of the Banach space$X$ and that of its dual $X^\ast$. Furthermore, we denote the space of bounded operators between normed spaces $X,Y$ by $B(X,Y)$ and write $\|\cdot\|_{X\to Y}$ for the operator norm.

\subsection{Tensor products and Hilbert-Schmidt operators}

Tensor products and simple tensors are the essential objects to understand kernel theorems for operators. See \cite{Ryan_2002} for an overview of the subject. In this paper, we use the following convention to define the \emph{simple tensor} of $f_1\in \H_1$ and $f_2\in \H_2$: 
$f_1 \otimes f_2$ will be understood as the rank one  operator mapping $\H_1$ on $\H_2$ according to
\begin{equation}
(f_1  \otimes f_2) (f) := \langle f ,  f_1 \rangle_{\H_1}  f_2  .
\label{eq:Tensorprodukt_als_Operator}
\end{equation}
Note that this tensor is homogeneous in the following sense: $\alpha (f_1\otimes f_2)=(\overline{\alpha}f_1)\otimes f_2=f_1\otimes (\alpha f_2).$
The tensor product $\H_1\otimes \H_2$ is then defined as the completion of the linear span of all simple tensors with respect to the metric induced by the inner product 
\begin{equation}\label{inner-prod-simple-tensor}
\langle f_1\otimes f_2 ,  g_1\otimes g_2\rangle_{\H_1\otimes \H_2} =\overline{\langle f_1  ,  g_1\rangle_{\H_1}}\langle f_2  ,  g_2\rangle_{\H_2}. 
\end{equation}
 We note here that  the tensor product 
$\H_1\otimes \H_2$ can be identified with the space of \emph{Hilbert-Schmidt operators} $\HS(\H_1,\H_2)$, which is defined as the space of bounded linear operators mapping $\H_1$ on $\H_2$ and equipped with the inner product
\begin{equation*}
\langle O_1  , O_2 \rangle_{\HS(\H_1,\H_2)} := \sum_{i \in \mathbb{N}} \langle O_1 e_i  ,   O_2 e_i \rangle_{\H_2}
\label{eq:Hilbert_Schmidt-Skalarprodukt}
\end{equation*} 
where $\{ e_i \}_{i \in \mathbb{N}}$ is any orthonormal basis of the Hilbert space $\H_1$. We also note that this identification constitutes a non-trivial kernel theorem \cite{conway}. If $O\in \HS(\H_1,\H_2)$, $f_1\in \H_1$ and $f_2\in \H_2$, then 
\begin{equation*}
\langle O ,  f_1 \otimes f_2 \rangle_{\HS(\H_1,\H_2)} = \langle O f_1 ,   f_2 \rangle_{\H_2} .
\label{eq:Wirkung_Tensorprodukt_als_Funktional_auf_Operator}
\end{equation*}

\subsection{Frame theory}

\noindent Throughout this work we use  \emph{frames} to decompose or synthesise functions and operators. See \cite{Christensen_2008, Casazza_2013} for an overview of frame theory. Let us recollect the definition of a frame originally given in  \cite{Duffin_1952}.% and Banach spaces \cite{Groechenig_1991}.

\begin{definition} %Let $H$ be a Hilbert space.
A countable set $\Psi : = \{ \psi_i \}_{i \in I}\subset \H$   is called a \emph{frame} for $\H$ if there are   positive numbers $A_{\Psi}$ and $B_{\Psi}$, called \emph{lower} and \emph{upper frame bounds}, such that, for any $f\in \H$
\begin{equation}
A_{\Psi} \, \| f \|_\H^2 \leqslant
  \sum_{i \in I} \vert \langle f  ,  \psi_i \rangle_\H \vert^2 \leqslant
  B_{\Psi} \, \| f \|_\H^2.
\label{eq:Rahmenungleichung_Hilbert}
\end{equation} 
\label{eq:Rahmen_im_Hilbertraum}
\end{definition}
Let us list  a few implications of this definition \begin{comment} \eqref{eq:Rahmenungleichung_Hilbert}
\end{comment}
that we shall make use of later. The \emph{analysis operator} 
\begin{equation*}
C_{\Psi} : \H \rightarrow  \ell^2(I) ,\quad f \mapsto \{\langle f   ,   \psi_i \rangle_\H\}_{i \in I},
\label{eq:Analysenoperator}
\end{equation*} is bounded and injective and has a closed range. The   \emph{synthesis operator}
\begin{equation*}
D_\Psi : \ell^2(I) \rightarrow \H ,\quad \{c_i\}_{i \in I} \mapsto  \sum_{i \in I} c_i \psi_i,
\label{eq:Synthesenoperator}
\end{equation*}
is bounded and surjective. The \emph{Gram operator}
\begin{equation}
G_\Psi : \ell^2(I) \rightarrow \ell^2(I) ,\quad \{c_i\}_{i \in I} \mapsto  \sum_{i' \in I} \langle \psi_{i'} ,  \psi_i \rangle_\H  \, c_{i'}
= \sum_{i' \in I} G_{i, i'} \, c_{i'}
= C_{\Psi} D_{\Psi}c ,
\label{eq:Gramscher_Operator}
\end{equation} with the \emph{Gram matrix} given by
\begin{equation}
G_{\Psi} := (G_{i, i'})_{(i, i') \in I^2} := \big( \langle \psi_{i'} ,  \psi_{i} \rangle_\H \big)_{(i,  i') \in I^2}  ,
\label{eq:Gramsche_Matrix}
\end{equation} is bounded. Occasionally, we will also need the definition of the \emph{cross Gram matrix} 
\begin{equation}
G_{\Psi, \, \widetilde \Psi}  = \big( \langle \psi_{i'} , \, \widetilde{\psi}_i \rangle \big)_{(i, \, i') \in I^2}.
\label{eq:Kreuz-Gramsche_Matrix}
\end{equation}
There exists at least one other frame $\Psi^\mathbf{d}  : = \{ \psi_i^{d} \}_{i \in I}$ for $\H$, called a \emph{dual frame}, that, together with $\Psi$, satisfies
\begin{equation}
\begin{split}
f  
   = \sum_{i \in I} \langle f ,   \psi_i \rangle_\H
 \,   \psi_i^{d}
   = D_{ {\Psi}^\mathbf{d} } C_{\Psi} f 
   = \sum_{i \in I}  \langle f  ,  \psi_i^{d}  \rangle_\H \, \psi_i 
   = D_{\Psi} C_{ {\Psi^\mathbf{d} }} f.
\end{split}
\label{eq:Rahmen_im_Hilbert_Wiederaufbau}
\end{equation}
Finally the \emph{frame operator}
\begin{equation*}
S_\Psi : \H \rightarrow \H ,\quad f \mapsto  \sum_{i \in I} \langle f   ,   \psi_i \rangle_\H\, \psi_i,
\label{eq:Rahmenoperator}
\end{equation*} is bounded, self-adjoint and invertible for any frame $\Psi$ and the sequence $\widetilde{\Psi}=\{ \widetilde{\psi}_i \}_{i \in I}=\{ S_{\Psi}^{-1} \psi_i \}_{i \in I}$, known as the \emph{canonical dual frame}, is a dual frame for $\Psi$.

\subsection{Localised frames and their co-orbit spaces}

One of the major advantages of frames over bases is that the former allow much more freedom in their construction and so can often be designed in a way that allows sparse decompositions of the functions of a given class.  Any element of a Hilbert space can be expressed as a linear combination of the elements of a frame, see \eqref{eq:Rahmen_im_Hilbert_Wiederaufbau}. To guarantee this reconstruction property in a Banach space, a more refined notion of frames has to be used, namely that of a \emph{Banach frame} \cite{Groechenig_1991}. In order to construct a Banach frame that achieves a sparse decomposition of a concrete Banach space, say, a specific decomposition space \cite{Feichtinger_1985}, the theory of structured Banach frames \cite{Borup_2007, Nielsen_2012, Nielsen_2014, Voigtlaender _2022} can be employed and indeed proved very useful \cite{Bytchenkoff_2020, Bytchenkoff_2021}. The Banach frames of choice for  more abstract co-orbit  spaces are  the so-called \emph{localised frames}.

Let us start with an intuitive notion of localisation \cite{Groechenig_2003} that should serve as an inspiration for the general concept of $\mathcal{A}$-localisation. For a given a metric $\rho$ on the indexing set $I$ that satisfies
\begin{equation*}
\inf_{i, \, i' \in I; \, i \neq i'} \rho ( i, i') = C > 0,
\label{eq:getrennte_Index-Menge}
\end{equation*}
the frame is called  \emph{localised}  %set $\Psi = \{ \psi_i \}_{i \in I}$ of a Hilbert space $H$#
if the absolute values of the elements of the   { Gram matrix} $G_\Psi $  decay polynomially as they deviate from the main diagonal of the matrix \cite{Groechenig_2003}, i.e., if
\begin{equation}
|(G_\Psi)_{i,i'}\vert=\vert \langle \psi_i , \psi_{i'} \rangle_\H \vert \lesssim (1+ \rho (i, i'))^{-n},\quad i,i'\in I,
\label{eq:Localisation_1}
\end{equation} for  some $n > n_0$.
The matrices whose elements satisfy \eqref{eq:Localisation_1} belong to a spectral matrix algebra. This allows the generalization of the notion of localisation \cite{Groechenig_2004}.

\begin{definition}\label{def:spectral}
An involutive Banach algebra $\mathcal{A}$ of infinite matrices with norm $\|\cdot\|_{\mathcal{A}}$ is called a \emph{spectral matrix algebra} if
\begin{enumerate}
\item [(i)]
every $A \in \mathcal{A}$ defines a bounded operator on $\ell^2 (I)$, i.e. $\mathcal{A} \subset B (\ell^2 (I))$, 
\item[(ii)] $\mathcal{A}$ is inverse-closed  in $B (\ell^2 (I))$, i.e. if $A \in \mathcal{A}$ is invertible  on $B(\ell^2(I))$, then $A^{-1} \in \mathcal{A}$, and
\item[(iii)] $\mathcal{A}$ is $solid$, i.e. if $A \in \mathcal{A}$ and $\vert  b_{i,j} \vert \leqslant \vert a_{i,j} \vert$ for any $i,j \in I$, then $B \in \mathcal{A}$ and $\left\| B \right\|_{\mathcal{A}} \leqslant \left\| A \right\|_{\mathcal{A}}$. 
\end{enumerate}
\label{eq:Spektrale-Matrix-Algebra}
\end{definition}
Examples of solid spectral matrix algebras include, among others, the Jaffard class \cite{Jaffard_1990} (which corresponds to \eqref{eq:Localisation_1} if $\rho$ is the standard Euclidean metric), the Schur class \cite{Schur} and the Sj\"ostrand class \cite{Sjöstrand}.
\begin{definition}
Let $I  $ be a  countable index set. A frame $\Psi = \{ \psi_i \}_{i \in I}\subset \H$   is said to be $\mathcal{A}$-\emph{intrinsically localised}, or $\mathcal{A}$-\emph{self-localised}, or simply $\mathcal{A}$-\emph{localised}, if its Gram matrix $G_\Psi :=(\langle \psi_i , \psi_{i'} \rangle_\H)_{(i, ,i') \in I^2}$ belongs to a spectral matrix algebra $\mathcal{A}$. 
\label{eq:Localisation_2_}
\end{definition}

Having a \lf frame $\Psi$ for a Hilbert space $\H$ allows to build a full range of so-called \emph{co-orbit spaces}. To do so, we consider a weight $w=\{w_i\}_{i\in I},\ w_i>0$, and define the weighted $\ell^p_w(I)$-spaces as the space of sequences on $I$ for which $\|c\|_{\ell^p_w(I)}=\|c\cdot w\|_{\ell^p(I)}<\infty,\ 1\leq p\leq\infty.$  

\begin{definition}
Let $\mathcal{A}$ be a spectral matrix algebra indexed by   $I \times I$. A weight sequence $w := \{ w_i \}_{i \in I}\subset \R^+$ is called $\mathcal{A}$-\emph{admissible} if every $A \in \mathcal{A}$ defines a bounded operator on $\ell_w^p(I)$  for every $1 \leqslant p \leqslant \infty$.
%, i.e.
%\begin{equation*}
%\mathcal{A} 
%   \subseteq \bigcap_{1 \leqslant p \leqslant \infty} B \left( \ell_w^p (I)\right).\label{eq:Zugelassenes_Gewicht}
%\end{equation*}
%\label{eq:Zugelassenes_Gewicht_}
\end{definition}

\begin{definition}
Let $\mathcal{A}$ be a spectral matrix algebra,   $\Psi : = \{ \psi_i \}_{i \in I}\subset\H$ be an $\mathcal{A}$-localised frame, $\widetilde\Psi : = \{ \widetilde\psi_i \}_{i \in I}\subset\H$ its canonical dual frame, $w := \{ w_i \}_{i \in I}$  an $\mathcal{A}$-admissible weight and
\begin{equation}
\begin{split}
\H_{00}(\Psi) &
   = \left\lbrace  \sum_{i \in I} c_i \, \psi_i :  \   \{ c_i \}_{i \in I} \in c_{00}  \right\rbrace,
\end{split}
\label{eq:H_00}
\end{equation} 
where $c_{00}$ denotes the space of complex sequences with finitely many non-zero terms. For $1\le p<\infty,$ the \emph{co-orbit space} $\H_w^p (\Psi)$ is defined as the norm completion of $\H_{00}(\Psi)$ with respect to the norm $\|f\|_{\H^p_w(\Psi)}:=\left\| C_{\widetilde\Psi} f   \right\|_{\ell_w^p}$,
%and  is called \emph{co-orbit space} generated by $\Psi$, 
while for $p=\infty$, $\H^\infty_w(\Psi)$ is defined as the completion of $\H_{00}(\Psi)$ with respect to the $\sigma(\H,\H_{00}(\Psi))$-topology. If $w\equiv 1$ we write $\H^p(\Psi)=\H^p_1(\Psi).$
\label{eq:Coorbit-Raum}
\end{definition}
\begin{remark}
More specifically, $\H^\infty_w(\Psi)$ is defined as a set of equivalence classes of sequences from $\H$, see \cite[Definition~3]{Balazs_2017}. We do not provide any details   for the moment, but shall do so when we consider co-orbit spaces generated by tensor products of \lf frames in Section~\ref{sec:tensor-of-frames}.

Also note that, while in \cite{Balazs_2017,Fornasier_2005} one may use any dual frame $\Psi^\mathbf{d}$ to define the co-orbit norms, we restrict ourselves to the canonical dual for the sake of simplicity.
\end{remark}

In the following proposition, we summarise the main properties of co-orbit spaces that we shall need. See \cite{Balazs_2017,Fornasier_2005} for the corresponding original statements. 
\begin{proposition}\label{prop:properties}
Let $\Psi$ and $\Phi$ be two $\mathcal{A}$-localised  frames for a Hilbert space $\H$,  $w$ be an $\mathcal{A}$-admissible weight, and $1\le p\le \infty$. Then the following statements hold:
\begin{enumerate}[label=(\roman*)]
%    \item\label{bp:1} the co-orbit spaces $\H_w^p (\Psi)$ are Banach spaces; (trivial)
    
    \item \label{bp:2}  $D_\Psi:\ell^p_w(I)\to \H^p_w(\Psi)$ is continuous;
    
    \item \label{bp:3}  $\H = \H_1^2(\Psi)$;
    
    \item\label{bp:4}  duality: for $1\le p<\infty,$   and $q$ satisfying $1/p+1/q=1$, one has $(\H_w^p(\Psi))^\ast=\H^q_{1/w}(\Psi)$,  where the duality is defined by 
$$
\langle f,g\rangle_{\H^p_w(\Psi),\H^q_{1/w}(\Psi)}=\langle C_{\widetilde{\Psi}}f,C_{{\Psi}}g\rangle_{\ell^p_w(I),\ell^q_{1/w}(I)};
$$ 
\item\label{bp:5} inclusions: $\H_{w_1}^{p_1} (\Psi) \subset \H_{w_2}^{p_2} (\Psi)$ for $w_1\geqslant w_2$, and $1 \leqslant p_1 < p_2 \leqslant \infty;$

\item \label{bp:7} $G_{\widetilde{\Psi}}\in\mathcal{A}$ and $G_{\widetilde{\Psi},\Psi}\in\mathcal{A}$, i.e.,
 the canonical dual frame $\widetilde{\Psi}$ is also $\mathcal{A}$-localised;
 
\item\label{bp:6} independence of $\Psi$:  if $G_{\Psi,\widetilde{\Phi}}\in\mathcal{A}$ and $G_{\widetilde{\Psi} ,\Phi}\in\mathcal{A}$, then 
 $\H^p_w(\Psi)=\H^p_w({\Phi})$ with equivalent norms; in particular, $\H^p_w(\Psi)=\H^p_w\big(\widetilde{\Psi}\big)$.
\end{enumerate}
\end{proposition}
The following lemma sets an upper bound of the $\H^1_w(\Psi)$-norm of the elements of $\Psi$ and $\widetilde{\Psi}$.

\begin{lemma}\label{lem:bound-norm-frame-elements}
Let $1\leq p<\infty$, $\Psi$ be an $\mathcal{A}$-localised frame for $\H$, and $w$ an $\mathcal{A}$-admissible weight. Then
\begin{equation}\label{eq:bound-norm-frame-elements}
    \|\psi_i\|_{\H^p_w(\Psi)}\lesssim w_i,\quad \text{and}\quad     \big\|\widetilde{\psi}_i\big\|_{\H^p_w(\Psi)}\lesssim w_i,\quad i\in I.
\end{equation}
\end{lemma}
\proof 
Remember that the Gram matrix $G_{\widetilde \Psi}$ of $\widetilde \Psi$ as well as the cross-Gram matrix of $\Psi$ and $\widetilde \Psi$ are bounded operators on $\ell_{w}^p(I)$ for $p \in [1, \, \infty]$, as  $\Psi$ and $\widetilde \Psi$ are  $\mathcal{A}$-localised by assumption and Proposition~\ref{prop:properties}~\ref{bp:7}, and $w$ is $\mathcal{A}$-admissible. Therefore,
\begin{equation*}
\begin{split}
 \big\| \widetilde \psi_{ i} \big\|_{\H_{w}^{p}\left( \Psi\right)}^p &
  = \sum_{k \in I} \big\vert \big\langle \widetilde \psi_{ i} ,\widetilde \psi_{ k} \big\rangle_\H \big\vert^p \, w_{ k} ^p 
  =  \sum_{k \in I} \big\vert \big( G_{\widetilde \Psi} \big)_{k,\, i} \big\vert^p \, w_{ k}^p
 \\ & =  \sum_{k \in I} \big\vert \big( G_{\widetilde \Psi} e_i \big)_k \big\vert^p \, w_{ k}^p
   = \left\| G_{\widetilde \Psi} e_i \right\|_{\ell_{w}^{p}(I)}^p
  \lesssim \left\| e_i \right\|_{\ell_{w}^{p}(I)}^p
=w_{i}^p,
 \end{split}
 \end{equation*} 
where $e_i$ stands for the $i$-th vector of the standard basis of $\ell^1_w(I)$. The upper bound of $\|  \psi_{ i} \|_{\H_{w}^{p}(\Psi)}$ follows from a similar argument, replacing $G_{\widetilde\Psi}$ by $G_{\widetilde{\Psi},\Psi}$.
\pbox

%\begin{equation}
%H_w^{p_1} (\Psi) \subset H_w^{p_2} (\Psi),\quad  1 \leqslant p_1 \leqslant p_2 \leqslant \infty,
%\label{eq:Coorbit-Raum-Skala}
%\end{equation}
%and $H_w^{p_1} (\Psi)$ is dense in $H_w^{p_2} (\Psi)$ \textcolor{blue}{as long as $1 \leqslant p_1 \leqslant p_2 < \infty$.} In other words, for any given weight $w$, the co-orbit spaces $H_w^{p} (\Psi)$ with $p \in    [1, \, \infty]$ constitute a scale of Banach spaces embedded in one another. Moreover, it holds \textcolor{blue}{$H = H_1^2(\Psi)$, and $(H_w^p(\Psi))^\ast=H^q_{1/w}(\Psi)$, $1/p+1/q=1$, $1\le p<\infty$. } Finally, if the frame $\Psi$   is $\mathcal{A}$-localised, then it constitutes a Banach frame for any element of the scale too \cite{Fornasier_2005}.  

 {Finally,  we recollect a result from  \cite{Fornasier_2005} concerning the representation of elements of $\H^1_w(\Psi)$ as a series expansion.}

\begin{proposition}\label{thm0} Let $\Psi$ be an $\mathcal{A}$-localised frame for $\H$ and $w$ be an $\mathcal{A}$-admissible weight.  For any $f \in \H_{w}^{1} \left( \Psi  \right)$ there is  a sequence of numbers $c := \left\lbrace c_{i} \right\rbrace_{i \in I} \in \ell_{w}^{1} (I)$ such that
\begin{equation*}
f = \sum_{i \in I} c_{i} \, \psi_{ i},
\end{equation*} where the series converges absolutely and
$
 \left\| c \right\|_{\ell_{w}^{1} \left( I \right)} \asymp  \left\| f \right\|_{\H_{w}^{1} \left( \Psi  \right)}.
$
\end{proposition}

\section{Tensor products of localised frames}\label{sec:tensor-of-frames}

If we consider the tensor product $\Psi_1\otimes \Psi_2=\{\psi_{1,\, i}\otimes\psi_{2,\, j}\}_{(i,j)\in I\times J}$ of two frames $\Psi_1:=\{\psi_{1,\, i}\}_{i\in I} \subset\H_1$ and $\Psi_2:=\{\psi_{2,\, j}\}_{j\in J} \subset\H_2$, then it is straightforward to show \cite{Balazs_2008_bis} that $\Psi_1\otimes \Psi_2$  forms a frame for the tensor product $\H_1\otimes \H_2$ and that the canonical dual frame is given by $\widetilde{\Psi_1\otimes\Psi_2}=\widetilde{\Psi}_1\otimes \widetilde{\Psi}_2$. If both frames are localised with respect to a solid spectral matrix algebra, then it is possible to also define co-orbit spaces. For the sake of completeness, we {provide here all necessary} definitions. Let us first define the mixed-norm sequence spaces
$$
\ell_w^{p,q}(I\times J):=\left\{c:I\times J\to \C:\ \sum_{j\in J}\Big(\sum_{i\in I} |c_{i,j}|^pw_{i,j}^p\Big)^{q/p}<\infty\right\} 
$$
for $1\le p,q\le \infty$, with the usual adaptations if $p=\infty,$ or $q=\infty.$ Since the order of summation matters if $p\neq q$, we also introduce the spaces
\begin{equation*}
\mathfrak{l}_{w}^{\, p, q} (I\times J)
   := \left\lbrace c: I \times J \to \mathbb{C}  : \quad \sum_{i \in I}
\Big( \sum_{j \in J}  \left\vert  c_{i, j} \right\vert^p  w_{i,  j}^p \Big)^{q/p} < \infty  \right\rbrace.
\label{eq:supergeile_Kernraeume_Spalte}
\end{equation*} 
Let
$$
\H_{00}(\Psi_1\otimes \Psi_2)=\left\{\sum_{(i,j)\in I\times J}c_{i,\, j}\, \psi_{1,\, i}\otimes\psi_{2,\, j}:\ \{c_{i,\, j}\}_{(i,j)\in I\times J}\in c_{00}\right\},
$$
and $p,q \in [1, \, \infty)$. The \emph{co-orbit spaces} $\H^{p, \, q}_w(\Psi_1\otimes \Psi_2)$ are defined as the completion of $\H_{00}(\Psi_1\otimes\Psi_2)$ with respect to the norm $\|F\|_{\H^{p,q}_w(\Psi_1\otimes \Psi_2)}:=\|C_{\widetilde{\Psi_1\otimes \Psi_2}}F\|_{\ell ^{p,q}_w(I\times J)}$, and  $\mathfrak{H}^{p, \, q}_w(\Psi_1\otimes \Psi_2)$ as the completion of $\H_{00}(\Psi_1\otimes\Psi_2)$ with respect to the norm $\|F\|_{\mathfrak{H}^{p,q}_w(\Psi_1\otimes \Psi_2)}:=\|C_{\widetilde{\Psi_1\otimes \Psi_2}}F\|_{\mathfrak{l} ^{p,q}_w(I\times J)}$. If $w\equiv 1$, we write $\H^{p,q}(\Psi_1\otimes\Psi_2)=\H^{p,q}_1(\Psi_1\otimes\Psi_2) $, and $\mathfrak{H}^{p,q}(\Psi_1\otimes\Psi_2)=\mathfrak{H}^{p,q}_1(\Psi_1\otimes\Psi_2) $. {We shall primarily consider the case where} $p=q$, but $p\neq q$ will make \textcolor{blue}{an} appearance in Theorem~\ref{Supergeiler_Satz}.

As mentioned in \cite{Fornasier_2005}, dealing with $\H_w^{\infty} (\Psi)$ for a localised frame $\Psi$ is technically rather tricky; a detailed exposition of the construction can be found in \cite{Balazs_2017}. Here we have to grapple with this issue again as now the tensor product frame cannot be assumed to be localised. The major obstacle to simply applying established theory to the tensor product of an $\mathcal{A}_1$-localised frame and an $\mathcal{A}_2$-localised frame is the fact that the natural candidate for a tensor product algebra, the projective tensor product $\mathcal{A}_1 \widehat\otimes_\pi \mathcal{A}_2$, 
%for the algebra to which the Gram matrix of the tensor product frame $\Psi_1 \otimes \Psi_2$ would have to belong 
is in general neither solid nor inverse-closed.

In order to properly define the co-orbit spaces where one of the indices or both of them are $\infty$, we need the following notion. We call 
two sequences $\left\lbrace F_{n } \right\rbrace_{ n \in \mathbb{N}}$, $\left\lbrace F'_{n} \right\rbrace_{n \in \mathbb{N}}\subset \H_1 \otimes \H_2$ \emph{ equivalent}  if
\begin{equation*}
\lim_{n  \to \infty} \big\langle F_{n } - F'_{n } , \widetilde \psi_{1, \, i} \otimes \widetilde \psi_{2, \, j} \big\rangle_{\H_1\otimes \H_2} = 0,
\label{eq:equi-Fol}
\end{equation*}
for all $(i, j) \in I \times J$.
\label{eq:equivalente-Folgen}

\begin{definition} Let $1\le p< \infty$, $n=1,  2,$ and $w_n$ be $\mathcal{A}_n$-admissible weights. We define the spaces $\H_{w_1 \otimes w_2}^{r,s} \left( \Psi_1 \otimes \Psi_2 \right)$ and $\mathfrak{H}_{w_1 \otimes w_2}^{r,s} \left( \Psi_1 \otimes \Psi_2 \right)$ where $(r,s)\in\{(p,\infty),\\ (\infty,p),(\infty,\infty)\}$, as the spaces of  equivalence classes   $F = \left[ \left\lbrace F_{n } \right\rbrace_{ n  \in \mathbb{N}} \right]$ with   $\left\lbrace F_{n } \right\rbrace_{ n  \in \mathbb{N} }\\ \subset  \H_1 \otimes \H_2$ such that  
\begin{equation}
\lim_{n \to \infty} \big\langle F_{n},  \widetilde \psi_{1, \, i} \otimes \widetilde \psi_{2, \, j} \big\rangle_{\H_1\otimes \H_2} = : k_{i,\, j}
\label{eq:Beingung_1}
\end{equation}
exists for all $(i,  j) \in I \times J$ and
\begin{equation*}
\sup_{n } \big\| C_{\widetilde{\Psi_1 \otimes \Psi_2}}   F_{n } \big\|_{\ell_{w_1 \otimes w_2}^{r, s} \left( I \times J \right)} < \infty,\quad \text{and}\quad\sup_{n } \big\| C_{\widetilde{\Psi_1 \otimes \Psi_2}}   F_{n } \big\|_{\mathfrak{l}_{w_1 \otimes w_2}^{r, s} \left( I \times J \right)} < \infty
\label{eq:Beingung_2}
\end{equation*}
respectively.
Furthermore we {define} $\big(C_{\widetilde{\Psi_1\otimes\Psi_2}}F\big)_{i,j}:=k_{i,j}$,  $(i,  j) \in I \times J$,
and
\begin{equation}
\left\| F \right\|_{\H_{w_1 \otimes w_2}^{r, s} \left( \Psi_1 \otimes \Psi_2 \right)}
 := \big\| C_{\widetilde{\Psi_1\otimes\Psi_2}}F\big\|_{\ell_{w_1 \otimes w_2}^{r,  s} \left( I \times J \right)}
,
\label{eq:Norm_Hinf-inf_w-w}
\end{equation}
as well as 
\begin{equation}
\left\| F \right\|_{\mathfrak{H}_{w_1 \otimes w_2}^{r, s} \left( \Psi_1 \otimes \Psi_2 \right)}
 := \big\| C_{\widetilde{\Psi_1\otimes\Psi_2}}F\big\|_{\mathfrak{l}_{w_1 \otimes w_2}^{r,  s} \left( I \times J \right)}
.
\label{eq:Norm_Hinf-inf_w-w-2}
\end{equation}
% where $k := \left\lbrace k_{i,\, j} \right\rbrace_{(i, \, j) \in I \times J}$.
%\label{eq:Hinf-inf_w-w}
\end{definition}

%The limits, which define $k$ via (\ref{eq:Beingung_1}), are linear functionals while $\left\| \, \cdot \, \right\|_{\ell_{w_1 \otimes w_2}^{\infty,  \infty} \left( I \times J \right)}$, being a norm, is absolutely homogeneous and satisfies the triangle inequality. 
It is clear that $\left\| \, \cdot \, \right\|_{\H_{w_1 \otimes w_2}^{r,s} \left( \Psi_1 \otimes \Psi_2 \right)}$ is a seminorm. Moreover, $\left\| F \right\|_{\H_{w_1 \otimes w_2}^{r, s} \left( \Psi_1 \otimes \Psi_2 \right)} =0$ implies that all limits in \eqref{eq:Beingung_1} equal zero for any sequence $\left\lbrace F_{n} \right\rbrace_{n \in \mathbb{N}}$ representing $F$. In other words, all such sequences $\left\lbrace F_{n} \right\rbrace_{ n  \in \mathbb{N}}$ are equivalent to the sequence of zeros in $\H_1 \otimes \H_2$ and so $F=0$. Consequently, $\left\| \, \cdot \, \right\|_{\H_{w_1 \otimes w_2}^{r, s} \left( \Psi_1 \otimes \Psi_2 \right)}$ is indeed a norm. The same reasoning shows that $\|\cdot\|_{\mathfrak{H}^{r,s}_{w_1\otimes w_2}(\Psi_1\otimes\Psi_2)}$ is a norm.

%The condition (\ref{eq:Beingung_2}) implies that there is  a non-negative finite $C$ such that
%\begin{equation*}
%\begin{split}
%\big\vert \big\langle K_{n, \, m},  \widetilde \psi_{2, \, i} \otimes \widetilde \psi_{1, \, j} \big\rangle \big\vert   \, w_{2, \,i} \, w_{1, \, j} \leqslant C
 %\end{split}
%\label{eq:Koeffizienten_Erklaerung_0}
%\end{equation*} for any $(n, \, m) \in \mathbb{N}^2$ and any $(i, \, j) \in I \times J$. Therefore
%\begin{equation*}
%\vert k_{i, \, j} \vert \, w_{2, \,i} \, w_{1, \, j} = \lim_{n, \, m \to \infty} \big\vert \big\langle K_{n, \, m},  \widetilde \psi_{2, \, i} \otimes \widetilde \psi_{1, \, j} \big\rangle \big\vert  \, w_{2, \,i} \, w_{1, \, j} \leqslant C,
%\label{eq:Koeffizienten_Erklaerung}
%\end{equation*}
%and so $k \in \ell_{w_1 \otimes w_2}^{\infty, \infty} \left( I \times J \right)$. In other words $ C_{\widetilde{\Psi_1 \otimes \Psi_2}}$ is a bounded, indeed isometric, operator mapping $H_{w_1 \otimes w_2}^{\infty,   \infty} \left( \Psi_1 \otimes \Psi_2 \right)$ on $\ell_{w_1 \otimes w_2}^{\infty,  \infty} \left( I \times J \right)$.
%\newline

The  spaces $\H^{p,p}_{w_1,w_2}(\Psi_1\otimes\Psi_2)$ are structurally very similar  to the co-orbit spaces of Section~\ref{sec:background} and results analoguous to Propositions~\ref{prop:properties} and  \ref{thm0} hold. These results however are not a direct consequence of the existing theory as, e.g., Propositions~\ref{prop:properties} and  \ref{thm0} assume the localisation of the frame with respect to a solid spectral algebra. We conjecture that $\Psi_1 \otimes \Psi_2$ would be $\mathcal{A}_1 \widetilde{\otimes} \mathcal{A}_2$-localised (for an appropriately chosen  tensor product structure satisfying the assumptions of Definition~\ref{def:spectral}) given that $\Psi_1$ and $\Psi_2$ are $\mathcal{A}_1$- and $\mathcal{A}_2$-localised respectively. This would allow us to directly apply existing theory. This is however not a simple matter and is the subject of an ongoing project of ours. Fortunately, the assumption that {both} $\Psi_1$ and $\Psi_2$ are \lf {suffices} to prove all the statements {that we need here}.

\begin{proposition}\label{Zerlegungssatz_Satz}
Let $\Psi_n,\Phi_n\subset\H_n$ be  $\mathcal{A}_n$-localised frames, and  $w_n$  be $\mathcal{A}_n$-admissible weights,  $n=1,  2$,  and $1\le p\le \infty$. Then the following statements hold:
 \begin{enumerate}[label=(\roman*)]
 \item\label{en:i} if $G_{\Psi_n,\Phi_n}\in\mathcal{A}_n$, $n=1,2$, then the cross Gram matrix %$G_{\widetilde {\Psi_1 \otimes \Psi_2},   \Psi_1 \otimes \Psi_2} $,
$  G_{\Psi_1 \otimes \Psi_2,   \Phi_1 \otimes \Phi_2} 
$
defines a bounded operator on $\ell^{p,p}_{w_1\otimes w_2}(I\times J)$;

 \item\label{en:ii} independence of $\Psi_1\otimes \Psi_2$: if $G_{\Psi_n,\widetilde{\Phi}_n}\in\mathcal{A}_n$ and $G_{  {\Phi}_n,\widetilde{\Psi}_n}\in\mathcal{A}_n$,  $n=1,2,$ then $$\H_{w_1\otimes w_2}^{p,p}(\Psi_1\otimes\Psi_2)=\H_{w_1\otimes w_2}^{p,p}\big({\Phi}_1\otimes {\Phi}_2\big)$$
    with equivalent norms. In particular, $$\H_{w_1\otimes w_2}^{p,p}(\Psi_1\otimes\Psi_2)=\H_{w_1\otimes w_2}^{p,p}\big(\widetilde{\Psi}_1\otimes\widetilde{\Psi}_2\big);$$

     \item \label{prop:(i)} $D_{\Psi_1\otimes\Psi_2}:\ell^{p,p}_{w_1\otimes w_2}(I\times J)\to \H_{w_1\otimes w_2}^{p,p}(\Psi_1\otimes\Psi_2)$ is bounded and the expansion converges absolutely for $p<\infty$, and weak-$\ast$ unconditionally for $p=\infty$;

 \item \label{prop:(ii)} for any $F \in \H_{w_1 \otimes w_2}^{1,  1} \left( \Psi_1 \otimes \Psi_2 \right)$, there exists a sequence $c      \in \ell_{w_1 \otimes w_2}^{1,1} (I \times J)$ such that
$F =  D_{\Psi_1 \otimes \Psi_2}   c$ and $\| F \|_{\H_{w_1 \otimes w_2}^{1,  1} \left( \Psi_1 \otimes \Psi_2 \right)} \asymp \| c \|_{\ell_{w_1 \otimes w_2}^{1,1} (I \times J)}$;

    \item\label{en:dual} duality:
for $1\leq p<\infty$ and $q$ satisfying $1/{p} + 1/{q} =1$,  
\begin{equation*}
   \left( \H_{w_1 \otimes w_2}^{p,  p} (\Psi_1 \otimes \Psi_2) \right)^* 
   \cong \H_{1/ (w_1 \otimes w_2)}^{q, q} (\Psi_1 \otimes \Psi_2)
\end{equation*} with the scalar product defined by
\begin{align*}
    \langle F, G \rangle_{ \H_{w_1 \otimes w_2}^{p,   p} (\Psi_1 \otimes \Psi_2)   ,   \H_{1/ (w_1 \otimes w_2)}^{q,   q} (\Psi_1 \otimes \Psi_2) }&
   \\  &\hspace{-3cm}:= \langle C_{\widetilde {\Psi_1 \otimes \Psi_2}} F , C_{\Psi_1 \otimes \Psi_2} G \rangle_{\ell_{w_1 \otimes w_2}^{p,  p} (I \times J)   ,   \ell_{1/(w_1 \otimes w_2)}^{q,  q} (I \times J)} .
\end{align*}
  %  \item \label{prop:(ibis)} $C_{\Psi_1\otimes\Psi_2}=D^*_{\Psi_1\otimes\Psi_2}:  \H_{w_1\otimes w_2}^{p,p   }(\Psi_1\otimes\Psi_2) \to \ell^{p,p}_{w_1\otimes w_2}(I\times J)$ is bounded;
\end{enumerate}
\end{proposition}

\proof 
Ad \ref{en:i}: As the weight functions $w_n$ are $\mathcal{A}_n$-admissible, it follows that the cross Gram matrix $G_{ {\Psi}_n,\Phi_n}$ defines a bounded operator on $\ell^p_{w_n}(I_n)$ \cite{Balazs_2017}. Therefore, for $1\le p< \infty,$ and $c  \in \ell_{w_1 \otimes w_2}^{p,p} (I \times J)$, we obtain from \eqref{inner-prod-simple-tensor} that
\begin{align*}
\big\|G_{ {\Psi_1 \otimes \Psi_2},   \Phi_1 \otimes \Phi_2} c\big\|_{\ell^{p,p}_{w_1\otimes w_2}(I\times J)}^{p}
&=\sum_{j\in J}\sum_{i\in I}\hspace{-2pt}\Big|\hspace{-2pt}\sum_{(i',j')\in I\times J}\hspace{-2pt}\overline{(G_{ {\Psi}_1,\Phi_1})_{i,i'}}(G_{ {\Psi}_2,\Phi_2})_{j,j'} c_{i',j'}\Big|^{p} w_{1,  i}^{p}  w_{2,  j}^{p}
\\
& \lesssim \sum_{j\in J}\sum_{i\in I}\Big|\sum_{j'\in   J}  (G_{ {\Psi}_2,\Phi_2})_{j,j'}c_{i,j'}\Big|^{p} w_{1,\, i}^{p}  w_{2,\, j}^{p}\\
& \lesssim \sum_{j\in J}\sum_{i\in I}|c_{i,j} |^{p} w_{1,\, i}^{p}  w_{2,\, j}^{p}= \| c\|_{\ell^{p,p}_{w_1\otimes w_2}(I\times J)}^p.
\end{align*}
A similar argument  {can be used to prove} the statement  {in the case where} $p=\infty$.  {This} concludes  {the proof of} \ref{en:i}.

Ad \ref{en:ii}: Since $F=D_{\Psi_1\otimes\Psi_2}C_{\widetilde{\Psi_1\otimes\Psi_2}}F$ for every $F\in\H_{00}(\Psi_1\otimes\Psi_2)$, we see that the frame coefficients with respect to $\widetilde{\Phi}_1\otimes\widetilde{\Phi_2}$ are given by 
$$ C_{\widetilde{\Phi_1\otimes\Phi_2}}F=C_{\widetilde{\Phi_1\otimes\Phi_2}}D_{\Psi_1\otimes\Psi_2}C_{\widetilde{\Psi_1\otimes\Psi_2}}F=G_{\Psi_1\otimes\Psi_2,\widetilde{\Phi_1\otimes\Phi_2}} C_{\widetilde{\Psi_1\otimes\Psi_2}}F. $$
Using \ref{en:i} we then conclude that for every $F\in\H_{00}(\Psi_1\otimes\Psi_2)$
\begin{align*}
\|F\|_{\H^{p,p}_{w_1\otimes w_2}(\Phi_1\otimes \Phi_2)}&=\|C_{\widetilde{\Phi_1\otimes\Phi_2}}F\|_{\ell^{p,p}_{w_1\otimes w_2}(I\times J)}\lesssim \|C_{\widetilde{\Psi_1\otimes\Psi_2}}F\|_{\ell^{p,p}_{w_1\otimes w_2}(I\times J)}
\\
&=\|F\|_{\H^{p,p}_{w_1\otimes w_2}(\Psi_1\otimes \Psi_2)}.
\end{align*}
Therefore, as $\H_{00}(\Psi_1\otimes\Psi_2)$ is dense in $\H^{p,p}_{w_1\otimes w_2}(\Psi_1\otimes \Psi_2)$ it follows that  $\H^{p,p}_{w_1\otimes w_2}(\Psi_1\otimes \Psi_2)\subseteq \H^{p,p}_{w_1\otimes w_2}(\Phi_1\otimes \Phi_2)$. Exchanging the roles of $\Psi_1\otimes \Psi_2$ and $\Phi_1\otimes \Phi_2$ results in $\H^{p,p}_{w_1\otimes w_2}(\Phi_1\otimes \Phi_2)\subseteq \H^{p,p}_{w_1\otimes w_2}(\Psi_1\otimes \Psi_2)$ and the norm equivalence.

Ad \ref{prop:(i)}: The proofs of \ref{prop:(i)} and \ref{prop:(ii)} follow  closely  the  arguments of   \cite[Proposition~2.4]{Fornasier_2005} and \cite[Lemmas~2 and~7]{Balazs_2017} respectively.
First, let $c \in \ell_{w_1 \otimes w_2}^{p,p} (I \times J)$ and set 
$
F  = D_{\Psi_1\otimes\Psi_2} c.
$ 
%where $I'$ and $J'$ stand for finite subsets of  $I$ and $J$ respectively. Then $F \in \H_{00} \left( \Psi_1 \otimes \Psi_2 \right)$ and
%\begin{equation*}
%C_{\widetilde {\Psi_1 \otimes \Psi_2}} F  = C_{\widetilde {\Psi_1 \otimes \Psi_2}} D_{\Psi_1 \otimes \Psi_2} \,  c \big|_{I' \times J'} = G_{\widetilde {\Psi_1 \otimes \Psi_2},  \Psi_1 \otimes \Psi_2} \, c \big|_{I' \times J'}.
%\label{eq:Zerlegungssatz_toto2}
%\end{equation*}
%where  the entries of the cross-Gram matrix $G_{\widetilde{\Psi_1 \otimes \Psi_2}, \, \Psi_1 \otimes \Psi_2} $ are given by  
%\begin{align*}
% \big( G_{\widetilde{\Psi_1 \otimes \Psi_2},   \Psi_1 \otimes \Psi_2} \big)_{(i,   j, j' ,  i') \in I \times J \times J \times I}
 %  :&= \big\langle \widetilde{\psi_{1,\, i}\otimes \psi_{2,\, j}},\psi_{1,\, i'}\otimes \psi_{2,\, j'}\big\rangle 
 %  \\ 
 %  &=\overline{\big\langle   \widetilde \psi_{1, \, i}, \psi_{1, \, i'}  \big\rangle }\big\langle \widetilde \psi_{2, \, j}, \psi_{2, \, j'}  \big\rangle.
%\end{align*}
Then
\begin{align} \label{eq:Zerlegungssatz_toto3_bis}
\left\| F  \right\|_{\H_{w_1 \otimes w_2}^{p,  p} \left( \Psi_1 \otimes \Psi_2 \right)} &
 = \big\| C_{\widetilde{\Psi_1 \otimes \Psi_2}} F \big\|_{\ell_{w_1 \otimes w_2}^{p,  p} \left( I \times J \right)} \\ & = \big\| G_{\widetilde{\Psi_1 \otimes \Psi_2},   \Psi_1 \otimes \Psi_2} \, c   \big\|_{\ell_{w_1 \otimes w_2}^{p,p} \left( I \times J \right)}
 \lesssim \big\| c  \big\|_{\ell_{w_1 \otimes w_2}^{p,p} \left( I \times J \right)},
\notag
\end{align} 
where we used \ref{en:i}, and Proposition~\ref{prop:properties}~\ref{bp:7}. The absolute convergence for $p<\infty$ follows immediately from Lemma~\ref{lem:bound-norm-frame-elements}.

\begin{comment}To show the unconditional weak-$\ast$ convergence for $p=\infty$, take $\varepsilon>0$, $c\in\ell^{\infty,\infty}_{w_1\otimes w_2}(I\times J)$, and $H\in \H_{00}(\Psi_1\otimes \Psi_2)$. We choose a finite subset $V_0\subset I\times J$ such that $\sum_{(i,j)\notin V_0}|\big\langle H,\psi_{1,\, i}\otimes \psi_{2,\, j}\rangle_{\H_1\otimes \H_2} \big| (w_{1,\, i}\, w_{2,\, j})^{-1}<\varepsilon/\|c\|_{\ell^{\infty,\infty}_{w_1\otimes w_2}(I\times J)}$. For two finite  sets  $ V_1,V_2\subseteq I\times J$ with $V_0 \subseteq V_1$ and $V_0\subseteq V_2$, it holds $V_1\backslash V_2\cup V_2\backslash V_1 \subseteq (I\times J)\backslash V_0$, and consequently 
\begin{align*}
    \Big| \Big\langle \sum_{(i,j)\in V_1} c_{i,j}\psi_{1,\, i}&\otimes \psi_{2,\, j}- \sum_{(i,j)\in V_2} c_{i,j}\psi_{1,\, i}\otimes \psi_{2,\, j},H\Big\rangle_{\H_1\otimes \H_2}\Big|
    \\
    &\leq    \Big| \Big\langle \sum_{(i,j)\in V_1\backslash V_2\cup V_2\backslash V_1} c_{i,j}\psi_{1,\, i}\otimes \psi_{2,\, j},H\Big\rangle_{\H_1\otimes \H_2}\Big|
    \\
    &\le\|c\|_{\ell^{\infty,\infty}_{w_1\otimes w_2}(I\times J)} \hspace{-1pt}\sum_{(i,j)\notin V_0}\hspace{-2pt}\big|\big\langle H,\psi_{1,\, i}\otimes \psi_{2,\, j}\rangle_{\H_1\otimes \H_2} \big| \frac{1}{w_{1,\, i}\, w_{2,\, j} }
   \hspace{-1pt}<\hspace{-1pt}\varepsilon.
\end{align*}
Therefore, the series $D_{\Psi_1\otimes \Psi_2}c$ converges weak-$\ast$ unconditionally 
which concludes the proof of \ref{prop:(i)}.
\end{comment}
The proof of the unconditional weak-$\ast$ convergence for $p=\infty$ is postponed  {until} the end of the proof  {of this proposition.}

%as the cross-Gram matrix
% is  a bounded operator on $\ell_{w_1 \otimes w_2}^{1,1} \left( I \times J \right)$. This relation extends to all $c\in\ell^{1,1}_{w_1\otimes w_2}(I\times J)$ by density, i.e.
%\begin{equation}
%\left\| F \right\|_{\H_{w_1 \otimes w_2}^{1,1} \left( \Psi_1 \otimes \Psi_2 \right)}
 %= \big\| C_{\widetilde{\Psi_1 \otimes \Psi_2}} F \big\|_{\ell_{w_1 \otimes w_2}^{1,1} \left( I \times J \right)}
%  \lesssim  \left\| c \right\|_{\ell_{w_1 \otimes w_2}^{1,1} \left( I \times J \right)},
%\label{eq:Zerlegungssatz_toto3_bis}
%\end{equation} where $F \in \H_{w_1 \otimes w_2}^{1,   1} \left( \Psi_1 \otimes \Psi_2 \right)$ as defined by \eqref{eq:Zerlegungssatz}. 

Ad \ref{prop:(ii)}: Let $F$ be an element of  $\H_{w_1 \otimes w_2}^{1,1} \left( \Psi_1 \otimes \Psi_2 \right)$. Then $F$ is the limit of a sequence $\{ F_k \}_{k \in \mathbb{N}} \subseteq \H_{00} \left( \Psi_1 \otimes \Psi_2 \right)$ with respect to the $\H_{w_1 \otimes w_2}^{1,1} \left( \Psi_1 \otimes \Psi_2 \right)$-norm. Being convergent,  $\{ F_k \}_{k \in \mathbb{N}}$ is a Cauchy sequence in $\H_{w_1 \otimes w_2}^{1,1} \left( \Psi_1 \otimes \Psi_2 \right)$. From this and the first equality in \eqref{eq:Zerlegungssatz_toto3_bis} we infer that $\{ C_{\widetilde{\Psi_1 \otimes \Psi_2}} F_k \}_{k \in \mathbb{N}}$ is a Cauchy sequence in $\ell_{w_1 \otimes w_2}^{1,1} \left( I \times J \right)$ converging to $c=C_{\widetilde{\Psi_1\otimes\Psi_2}}F \in \ell_{w_1 \otimes w_2}^{1,1} \left( I \times J \right)$. Let $F' := D_{\Psi_1 \otimes \Psi_2} \, c$. From \ref{prop:(i)} we infer that $F'\in \H_{w_1 \otimes w_2}^{1,1} \left( \Psi_1 \otimes \Psi_2 \right)$. Using the fact that$$
C_{\widetilde{\Psi_1\otimes\Psi_2}}H=C_{\widetilde{\Psi_1\otimes\Psi_2}}D_{\Psi_1\otimes\Psi_2}C_{\widetilde{\Psi_1\otimes\Psi_2}}H
$$ 
for any $ H \in\H_{00}(\Psi_1\otimes\Psi_2)$
and the boundedness of $D_{\Psi_1\otimes \Psi_2}$ {results in}
\begin{equation*}
\begin{split}
\lim_{k \to \infty} \big\| F'& - F_k \big\|_{\H_{w_1 \otimes w_2}^{1,1} \left( \Psi_1 \otimes \Psi_2 \right)} 
 = \lim_{k \to \infty} \big\| C_{\widetilde{\Psi_1 \otimes \Psi_2}} (F' - F_k) \big\|_{\ell_{w_1 \otimes w_2}^{1,1} \left( I \times J \right)} \\
 & = \lim_{k \to \infty} \big\| C_{\widetilde{\Psi_1\otimes\Psi_2}}D_{\Psi_1\otimes\Psi_2}C_{\widetilde{\Psi_1\otimes\Psi_2}}(F- F_k )\big\|_{\ell_{w_1 \otimes w_2}^{1,1} \left( I \times J \right)}
 \\ &\lesssim \lim_{k \to \infty} \big\| C_{\widetilde{\Psi_1\otimes\Psi_2}}(F- F_k )\big\|_{\ell_{w_1 \otimes w_2}^{1,1} \left( I \times J \right)}
 =0.
\end{split}
\label{eq:Zerlegungssatz_toto4}
\end{equation*}
Therefore $F' = F = D_{\Psi_1 \otimes \Psi_2}   c$ and we have proved that $D_{\Psi_1 \otimes \Psi_2}$ is surjective. Finally, $D_{\Psi_1 \otimes \Psi_2}$ is a bijective operator mapping $\ell_{w_1 \otimes w_2}^{1,1} \left( I \times J \right) / \text{ker} \left( D_{\Psi_1 \otimes \Psi_2} \right)$ \textcolor{blue}{on} $\H_{w_1 \otimes w_2}^{1,1} \left( \Psi_1 \otimes \Psi_2 \right)$\textcolor{blue}{,} which shows that there exists $c\in \ell^{1,1}_{w_1\otimes w_2}(I\times J)$ such that $D_{\Psi_1\otimes \Psi_2}c=F$ and 
\begin{equation}
\left\| c \right\|_{\ell_{w_1 \otimes w_2}^{1,1} \left( I \times J \right)}
  \lesssim \left\| F \right\|_{\H_{w_1 \otimes w_2}^{1,  1} \left( \Psi_1 \otimes \Psi_2 \right)},
\label{eq:Zerlegungssatz_toto5}
\end{equation}
according to the inverse mapping theorem. Combining \eqref{eq:Zerlegungssatz_toto3_bis} and \eqref{eq:Zerlegungssatz_toto5} results in the norm equivalence of \ref{prop:(ii)}.

Ad \ref{en:dual}:  We follow the line of arguments of  \cite[Proposition 2]{Balazs_2017}. Let us define a sesquilinear form $\H_{w_1 \otimes w_2}^{p,p} (\Psi_1 \otimes \Psi_2)\times \H_{1/ (w_1 \otimes w_2)}^{q, q} (\Psi_1 \otimes \Psi_2)\to \C$ by
 \begin{align*}
     \left( F, G \right) \mapsto  \langle C_{\widetilde {\Psi_1 \otimes \Psi_2}} F , C_{\Psi_1 \otimes \Psi_2} G \rangle_{\ell_{w_1 \otimes w_2}^{p,  p} (I \times J)   ,   \ell_{1/(w_1 \otimes w_2)}^{q,  q} (I \times J)}. 
 \end{align*} 
 Fixing $G \in \H_{1/ (w_1 \otimes w_2)}^{q, q} (\Psi_1 \otimes \Psi_2)$ defines a % on $\H_{w_1 \otimes w_2}^{p,  p} (\Psi_1 \otimes \Psi_2) \times \H_{1/ (w_1 \otimes w_2)}^{q,  q} (\Psi_1 \otimes \Psi_2)$ makes it into a 
 bounded linear functional  $ \mathcal{W}(G) :\H_{w_1 \otimes w_2}^{p,  p} (\Psi_1 \otimes \Psi_2) \to \C $ according to 
 \begin{equation*}
     \mathcal{W}(G) (F)=\langle C_{\widetilde {\Psi_1 \otimes \Psi_2}} F , C_{\Psi_1 \otimes \Psi_2} G \rangle_{\ell_{w_1 \otimes w_2}^{p,  p} (I \times J)   ,   \ell_{1/(w_1 \otimes w_2)}^{q,  q} (I \times J)}.
 \end{equation*}
 Indeed, by \ref{en:ii}
 \begin{equation*}
 \begin{split}
\vert   \mathcal{W}(G)   (F) \vert
 & \leqslant \| C_{\widetilde {\Psi_1 \otimes \Psi_2}} F \|_{\ell_{w_1 \otimes w_2}^{p,  p} (I \times J)} \, \| C_{\Psi_1 \otimes \Psi_2} G \|_{\ell_{1/(w_1 \otimes w_2)}^{q, q} (I \times J)}
 \\
 &   =  \| F \|_{\H_{w_1 \otimes w_2}^{p,  p} (\Psi_1 \otimes \Psi_2)} \, \|   G \|_{\H_{1/(w_1 \otimes w_2)}^{q, q} (\widetilde{\Psi}_1 \otimes \widetilde{\Psi}_2)}
  \\
 &   \lesssim  \| F \|_{\H_{w_1 \otimes w_2}^{p,  p} (\Psi_1 \otimes \Psi_2)} \, \|   G \|_{\H_{1/(w_1 \otimes w_2)}^{q, q} ( {\Psi_1} \otimes  {\Psi_2})},
     \end{split}
 \end{equation*}
 for any $F \in \H_{w_1 \otimes w_2}^{p,  p} (\Psi_1 \otimes \Psi_2)$. This shows  that      $\mathcal{W}:\H_{1/ (w_1 \otimes w_2)}^{q,  q} (\Psi_1 \otimes \Psi_2)\to  \left( \H_{w_1 \otimes w_2}^{p,  p} (\Psi_1 \otimes \Psi_2) \right)^*$ is bounded as
 \begin{equation*}
\|   \mathcal{W}(G) \|_{  \left( \H_{w_1 \otimes w_2}^{p,  p} (\Psi_1 \otimes \Psi_2) \right)^* }
  \lesssim \|  G \|_{\H_{1/(w_1 \otimes w_2)}^{q,  q} (\Psi_1\otimes \Psi_2)}.
 \end{equation*}

 Conversely, let $ H \in \left( \H_{w_1 \otimes w_2}^{p, p} (\Psi_1 \otimes \Psi_2) \right)^*$ and $c  \in \ell_{w_1 \otimes w_2}^{p,  p} (I \times J)$ with $p \in  [1, \infty )$. Then $D_{\Psi_1 \otimes \Psi_2} c \in \H_{w_1 \otimes w_2}^{p,  p} (\Psi_1 \otimes \Psi_2)$ by \ref{prop:(i)} and
 \begin{equation*}
 H \left( D_{\Psi_1 \otimes \Psi_2} c \right)
  = \sum_{i \in I}  \sum_{j \in J}  c_{i,  j} H (\psi_{1, \, i} \otimes \psi_{2, \, j})
 \end{equation*} where the series converges unconditionally for any $c \in \ell_{1/{w_1 \otimes w_2}}^{p, \, q} (I \times J)$.   Therefore $\{ H (\psi_{1, \, i} \otimes \psi_{2, \, j}) \}_{(i,  j) \in I \times J} \in \ell_{1/(w_1 \otimes w_2)}^{q, \, q} (I \times J)$ \cite{conway}. If we define the mapping 
 $
     \mathcal{V} : \left( \H_{w_1 \otimes w_2}^{p, \, p} (\Psi_1 \otimes \Psi_2) \right)^* \to \H_{1/ (w_1 \otimes w_2)}^{q,  q} (\Psi_1 \otimes \Psi_2)
$  by
 \begin{equation*}
     \mathcal{V}(H) := \sum_{i \in I}  \sum_{j \in J}
     \overline{H (\psi_{1, \, i} \otimes \psi_{2, \, j})} \,
     \widetilde{\psi}_{1, \, i} \otimes \widetilde{\psi}_{2, \, j},\quad  H \in \left( \H_{w_1 \otimes w_2}^{p,  p} (\Psi_1 \otimes \Psi_2) \right)^*,
 \end{equation*} 
 then, for  $F \in \H_{w_1 \otimes w_2}^{p,  p} (\Psi_1 \otimes \Psi_2)$,
 \begin{align*}
( \mathcal{W }&( \mathcal{V } ( H  )  )  ) (F)  
= \langle C_{\widetilde {\Psi_1 \otimes \Psi_2}} F , C_{\Psi_1 \otimes \Psi_2} \mathcal{V}  ( H  ) \rangle_{\ell_{w_1 \otimes w_2}^{p, p} (I \times J) ,  \ell_{1/(w_1 \otimes w_2)}^{q,  q} (I \times J)} \\
& = \sum_{i \in I}  \sum_{j \in J} \langle F, \widetilde{\psi}_{1, \, i} \otimes \widetilde{\psi}_{2, \, j} \rangle \overline{\Big\langle \sum_{k \in I}  \sum_{l \in J} \overline{H (\psi_{1, \, k} \otimes \psi_{2, \, l})} \, \widetilde{\psi}_{1, \, k} \otimes \widetilde{\psi}_{2, \, l}   ,   \psi_{1, \, i} \otimes \psi_{2, \, j} \Big\rangle} \\
& = \sum_{k \in I}  \sum_{l \in J} H (\psi_{1, \, k} \otimes \psi_{2, \, l})
\Big\langle \sum_{i \in I}  \sum_{j \in J} \langle F, \widetilde{\psi}_{1, \, i} \otimes \widetilde{\psi}_{2, \, j} \rangle \psi_{1, \, i} \otimes \psi_{2, \, j}   ,   \widetilde{\psi}_{1, \, k} \otimes \widetilde{\psi}_{2, \, l} \Big\rangle \\
& = \sum_{k \in I}  \sum_{l \in J} H (\psi_{1, \, k} \otimes \psi_{2, \, l}) \langle F, \widetilde{\psi}_{1, \, k} \otimes \widetilde{\psi}_{2, \, l} \rangle \\
& = H \left( \sum_{k \in I}  \sum_{l \in J} \langle F, \widetilde{\psi}_{1, \, k} \otimes \widetilde{\psi}_{2, \, l} \rangle  \psi_{1, \, k} \otimes \psi_{2, \, l} \right)
= H(F),
 \end{align*} 
 where the change of the order of summations, taking scalar products and application of  $H$ is justified by the unconditional convergence of the series involved, the continuity of $H$, and a density argument. The operator $\mathcal{W}$ is therefore surjective. Moreover, for $G\in \H_{1/ (w_1 \otimes w_2)}^{q,  q} (\Psi_1 \otimes \Psi_2)$,
 \begin{align*}
         \mathcal{V} (\mathcal{ W} ( G ) ) &
         = \sum_{i \in I}  \sum_{j \in J}
     \overline{   \mathcal{W}  ( G  )   (\psi_{1, \, i} \otimes \psi_{2, \, j})} \,
     \widetilde{\psi}_{1, \, i} \otimes \widetilde{\psi}_{2, \, j}  \\
    & = \sum_{i \in I}  \sum_{j \in J} \langle G  ,  \psi_{1, \, i} \otimes \psi_{2, \, j} \rangle \widetilde{\psi}_{1, \, i} \otimes \widetilde{\psi}_{2, \, j}
    = G,
 \end{align*} 
which shows that $\mathcal{W}$ is also injective and  thus invertible. 

It remains to prove the unconditional weak-$\ast$ convergence of $D_{\Psi_1\otimes\Psi_2}$ when $p=\infty$.  Let $c\in\ell^{\infty,\infty}_{w_1\otimes w_2}(I\times J)$ and $H\in \H_{00}(\Psi_1\otimes \Psi_2)$. Then using \ref{en:ii} and \ref{en:dual} {results in}
\begin{align*}
 \big| \langle D_{\Psi_1\otimes\Psi_2}c,&H\rangle_{\H^{\infty,\infty}_{w_1\otimes w_2}(\Psi_1\otimes\Psi_2),\H^{1,1}_{1/(w_1\otimes w_2)}(\Psi_1\otimes\Psi_2)}\big| \\ & \leqslant \sum_{(i,j)\in I\times J}|c_{i,j}|\big|\langle C_{\widetilde{\Psi_1\otimes\Psi_2}} \psi_{1,\, i}\otimes \psi_{2,\, j},C_{{\Psi_1\otimes\Psi_2}}H\rangle_{\ell^{\infty,\infty}_{w_1\otimes w_2}(I\times J),\ell^{1,1}_{1/(w_1\otimes w_2)}(I\times J)} \big|
 \\  & = \sum_{(i,j)\in I\times J}|c_{i,j}| \big|(C_{\Psi_1\otimes\Psi_2}H)_{i,j}\big|\lesssim\|c \|_{\ell^{\infty,\infty}_{w_1\otimes w_2}(I\times J)}\|H\|_{\H^{1,1}_{1/(w_1\otimes w_2)}(\Psi_1\otimes \Psi_2)}.
\end{align*}
 By density {of $\H_{00}(\Psi_1\otimes \Psi_2)$ in $\H^{1,1}_{1/(w_1\otimes w_2)}(\Psi_1\otimes \Psi_2)$}, the chain of {inequalities} {also} holds for any $H\in \H^{1,1}_{1/(w_1\otimes w_2)}(\Psi_1\otimes \Psi_2)$\textcolor{blue}{,} showing that $D_{\Psi_1\otimes \Psi_2}c$ converges weak-$\ast$ absolutely and in particular  weak-$\ast$ unconditionally.
\pbox

 \begin{remark}
 It is unclear to us whether the statements of Proposition~\ref{Zerlegungssatz_Satz} apply to the co-orbit spaces $\H_{w_1 \otimes w_2}^{p_1,  p_2} (\Psi_1 \otimes \Psi_2) $ if $p_1\neq p_2$. The main problem is that the argument used in proving \ref{en:i} does not apply to mixed $\ell^{p,q}$-spaces. In particular, the cross Gram matrix $G_{\Psi_1\otimes \Psi_2,\Phi_1\otimes\Phi_2}$ is not necessarily a bounded operator on the spaces $\ell^{p_1,p_1}_{w_1\otimes w_2}(I\times J)$ when $p_1\neq p_2$. The proofs of the statements \ref{en:ii}-\ref{en:dual} rely crucially on \ref{en:i}.   As a consequence of one of our kernel theorems, we shall however show in Corollary~\ref{cor} that some mixed-norm co-orbit spaces are indeed independent of the specific frames defining them as long as the individual cross Gram matrices belong to  $\mathcal{A}_n$, $n=1,2$.
 \end{remark}

In the following we establish the reconstruction formula on $\H_{1/(w_1 \otimes w_2)}^{\infty,   \infty} \left( \Psi_1 \otimes \Psi_2 \right)$ and a necessary and sufficient condition for $k \in  \ell_{1/(w_1 \otimes w_2)}^{\infty, \, \infty} (I \times J)$ to be the generalised wavelet transform of a vector $F \in \H_{1/(w_1 \otimes w_2)}^{\infty,   \infty} \left( \Psi_1 \otimes \Psi_2 \right)$ -- which, following the tradition  of \cite{Feichtinger_1988,  Fornasier_2005}, can be called the 'correspondence principle'.
\newline

\begin{lemma}\label{thm3} Let $\Psi_n,\Phi_n\subset\H_n$ be  $\mathcal{A}_n$-localised frames, and  $w_n$  be $\mathcal{A}_n$-admissible weights for  $n=1,  2$. If $F \in \H_{1/(w_1 \otimes w_2)}^{\infty,   \infty} \left( \Psi_1 \otimes \Psi_2 \right)$, then 
\begin{equation}\label{eq:reco}
F=D_{\Psi_1\otimes\Psi_2}C_{\widetilde{\Psi_1\otimes\Psi_2}}F=D_{\widetilde{\Psi_1\otimes\Psi_2}}C_{\Psi_1\otimes\Psi_2}F.
\end{equation}

Moreover, for any  $k \in  \ell_{1/(w_1 \otimes w_2)}^{\infty,   \infty} (I \times J)$ there exists $F \in \H_{1/(w_1 \otimes w_2)}^{\infty,   \infty} \left( \Psi_1 \otimes \Psi_2 \right)$ such that
\begin{equation}
\begin{split}
 k = C_{\widetilde{\Psi_1 \otimes \Psi_2}} F
 \end{split}
\label{eq:Correspondence_1}
\end{equation}
if and only if
\begin{equation}
 k = G_{\widetilde{\Psi_1 \otimes \Psi_2}, \, \Psi_1 \otimes \Psi_2} \, k = C_{\widetilde{\Psi_1 \otimes \Psi_2}} D_{\Psi_1 \otimes \Psi_2} \, k.
\label{eq:Correspondence_2}
\end{equation}
The statements also apply to $k \in \ell_{w_1 \otimes w_2}^{\infty,  \infty} \left( I \times J \right)$ and $F \in \H_{w_1 \otimes w_2}^{\infty,   \infty} \left( \Psi_1 \otimes \Psi_2 \right)$.
\end{lemma}

\proof Let  $F \in \H_{1/(w_1 \otimes w_2)}^{\infty,   \infty} \left( \Psi_1 \otimes \Psi_2 \right)$, $H\in \H_{00}(\Psi_1\otimes\Psi_2)$, and let $\left\lbrace I_n \right\rbrace_{n \in \mathbb{N}}$ and $\left\lbrace J_n \right\rbrace_{n \in \mathbb{N}}$ be  sequences of finite subsets of $I$ and $J$ respectively such that $ I_n \subset I_{n+1} \subseteq I$, and   $ J_n \subset J_{n+1} \subseteq J$ for any $n \in \mathbb{N}$, $\cup_{n \in \mathbb{N}} I_n = I$ and $\cup_{n \in \mathbb{N}} J_n = J$. By Proposition~\ref{Zerlegungssatz_Satz}~\ref{prop:(i)} $D_{\Psi_1\otimes\Psi_2}C_{\widetilde{\Psi_1\otimes\Psi_2}}F$ is well-defined and converges weak-$\ast$ unconditionally. Moreover, 
\begin{align*}
\langle &D_{ {\Psi_1\otimes\Psi_2}}C_{\widetilde{\Psi_1\otimes \Psi_2}}F,H\rangle_{\H_{1/(w_1 \otimes w_2)}^{\infty,   \infty} \left( \Psi_1 \otimes \Psi_2 \right),\H_{w_1 \otimes w_2}^{1,1} \left( \Psi_1 \otimes \Psi_2 \right)}
\\
& =\lim_{n\to \infty}  \sum_{(i,j)\in I_n\times J_n}\hspace{
-3pt}\langle F,\widetilde{\psi}_{1,\, i}\otimes \widetilde{\psi}_{2,\, j}\rangle_{\H_{1/(w_1 \otimes w_2)}^{\infty,   \infty} \left( \Psi_1 \otimes \Psi_2 \right),\H_{w_1 \otimes w_2}^{1,1} \left( \Psi_1 \otimes \Psi_2 \right)}\langle \psi_{1,\, i}\otimes\psi_{2,\, j}, H \rangle%_{\H_1\otimes \H_2}
\\
&=\big\langle F, D_{\widetilde{\Psi_1\otimes\Psi_2}}C_{\Psi_1\otimes \Psi_2}H\big\rangle_{\H_{1/(w_1 \otimes w_2)}^{\infty,   \infty} \left( \Psi_1 \otimes \Psi_2 \right),\H_{w_1 \otimes w_2}^{1,1} \left( \Psi_1 \otimes \Psi_2 \right)} 
\\
&=\langle F,H\rangle_{\H_{1/(w_1 \otimes w_2)}^{\infty,   \infty} \left( \Psi_1 \otimes \Psi_2 \right),\H_{w_1 \otimes w_2}^{1,1} \left( \Psi_1 \otimes \Psi_2 \right)}.
\end{align*}
And so $D_{ {\Psi_1\otimes\Psi_2}}C_{\widetilde{\Psi_1\otimes \Psi_2}}F=F$.
The second reconstruction formula follows from a similar argument once we recall that $\H_{1/(w_1 \otimes w_2)}^{\infty,   \infty} \left( \Psi_1 \otimes \Psi_2 \right)=\H_{1/(w_1 \otimes w_2)}^{\infty,   \infty} \big( \widetilde{\Psi}_1 \otimes \widetilde{\Psi}_2 \big)$.

To prove the second statement, let $k \in  \ell_{1/(w_1 \otimes w_2)}^{\infty,  \infty} (I \times J)$   satisfy (\ref{eq:Correspondence_2}). Then
\begin{equation*}
 F_{n } := \sum_{i \in I_n} \sum_{j \in J_n} k_{i,  j}\,   \psi_{1, \, i} \otimes \psi_{2, \, j} \in \H_1 \otimes \H_2,
\label{eq:Correspondence_Beweis_1}
\end{equation*} 
for any $n\in \mathbb{N}$ and
\begin{equation*}
\begin{split}
\lim_{n \to \infty} \big\langle F_{n},  \widetilde \psi_{1, \, r} \otimes \widetilde \psi_{2, \, s} \big\rangle_{\H_1\otimes\H_2} &
 = \lim_{n \to \infty} \Big\langle \sum_{i \in I_n} \sum_{j \in J_n} k_{i, \, j}\,  \psi_{1, \, i} \otimes \psi_{2, \, j},  \widetilde \psi_{1, \, r} \otimes \widetilde \psi_{2, \, s} \Big\rangle_{\H_1\otimes\H_2} \\ 
 & = \lim_{n \to \infty} \sum_{i \in I_n} \sum_{j \in J_n} k_{i, \, j} \, \big\langle \psi_{1, \, i} \otimes \psi_{2, \, j}, \, \widetilde \psi_{1, \, r} \otimes \widetilde \psi_{2, \, s} \big\rangle_{\H_1\otimes\H_2} \\
 & = \left( C_{\widetilde{\Psi_1 \otimes \Psi_2}}  D_{\Psi_1 \otimes \Psi_2} \, k \right)_{r, \, s} \\
 & = \left( G_{\widetilde{\Psi_1 \otimes \Psi_2}, \, \Psi_1 \otimes \Psi_2} \, k \right)_{r, \, s}
 = k_{r, \, s}.
 \end{split}
\label{eq:Correspondence_Beweis_2}
\end{equation*}
From this we infer that
\begin{equation*}
\big\langle F_{n},   \widetilde \psi_{1, \, r} \otimes \widetilde \psi_{2, \, s} \big\rangle_{\H_1\otimes\H_2} = \left( G_{\widetilde{\Psi_1 \otimes \Psi_2}, \, \Psi_1 \otimes \Psi_2} \, k \big|_{I_n \times J_n} \right)_{r, \, s}
\label{eq:Correspondence_Beweis_2_bla}
\end{equation*} 
and so
\begin{align*}
\sup_{r, \, s}  \big\vert \big\langle F_{n},  \widetilde \psi_{1, \, r} \otimes \widetilde \psi_{2, \, s} \big\rangle_{\H_1\otimes\H_2} \big\vert  (w_{1,\, r} \, w_{2,\, s})^{-1}
 &\lesssim  \big\| k \big|_{I_n \times J_n} \big\|_{\ell_{1/(w_1 \otimes w_2)}^{\infty,  \infty}(I\times J)}
 \\ &\leqslant \| k \|_{\ell_{1/(w_1 \otimes w_2) }^{\infty, \, \infty}(I\times J)}.
\label{eq:Correspondence_Beweis_2_bla-bla}
\end{align*}
Therefore  $F = \left[ F_{n} \right] \in \H_{1/(w_1 \otimes w_2)}^{\infty,   \infty} \left( \Psi_1 \otimes \Psi_2 \right)$ satisfies (\ref{eq:Correspondence_1}) and is, as an equivalence class, uniquely defined by $k$ via \eqref{eq:Correspondence_1}.

Conversely, let $k \in  \ell_{1/(w_1 \otimes w_2)}^{\infty,  \infty} (I \times J)$ and    $F \in \H_{ 1/(w_1 \otimes w_2)}^{\infty,  \infty} \left( \Psi_1 \otimes \Psi_2 \right)$ be such that    \eqref{eq:Correspondence_1} is satisfied. Proposition~\ref{Zerlegungssatz_Satz}~\ref{prop:(i)} indicates that $D_{\Psi_1 \otimes \Psi_2}$ is a bounded operator mapping $\ell_{1/(w_1 \otimes w_2)}^{\infty,   \infty} (I \times J)$ to $\H_{1/(w_1 \otimes w_2)}^{\infty,   \infty} \left( \Psi_1 \otimes \Psi_2 \right)$. Applying  $C_{\widetilde{\Psi_1\otimes \Psi_2}}D_{\Psi_1 \otimes \Psi_2}$ to both sides of \eqref{eq:Correspondence_1} and using \eqref{eq:reco} results in 
\begin{equation*}
C_{\widetilde{\Psi_1\otimes \Psi_2}}D_{\Psi_1 \otimes \Psi_2} k = C_{\widetilde{\Psi_1\otimes \Psi_2}}D_{\Psi_1 \otimes \Psi_2} C_{\widetilde{\Psi_1 \otimes \Psi_2}} F =C_{\widetilde{\Psi_1\otimes \Psi_2}}F=k.
\label{eq:Correspondence_1_modif}
\end{equation*} 
Finally, if $w_n$ is an $\mathcal{A}_n$-admissible weight then so is $1/w_n$, $n=1,2$. Therefore, the statement of the theorem also applies to $k \in \ell_{w_1 \otimes w_2}^{\infty,   \infty} \left( I \times J \right)$ and $F \in \H_{w_1 \otimes w_2}^{\infty,  \infty} \left( \Psi_1 \otimes \Psi_2 \right)$. \pbox
\newline

\section{Kernel Theorems}\label{sec:proofs}
\noindent We are now ready to state and prove our main results. To prove the outer kernel theorem, we essentially follow the proof strategy of \cite{Balazs_2019}, but use somewhat more involved arguments in those steps where we cannot rely on either the group structure or submultiplicativity of the weight function. The basic idea is to represent the operator by its Galerkin matrix and to derive the kernel of the operator from this representation.
\begin{theorem}\label{kernsatz} \emph{(Outer Kernel Theorem)} Let $\Psi_n\subset \H_n$ be an $\mathcal{A}_n$-localised frame, and $w_n$ an $\mathcal{A}_n$-admissible weight for $n=1,2$. 

To any element of $K \in \H_{1 / (w_1 \otimes w_2)}^{\infty,   \infty} \left( \Psi_1 \otimes \Psi_2 \right)$ corresponds a unique bounded linear operator $O:\H_{w_1}^{1} \left( \Psi_1 \right)\to \H_{1 / w_2}^{\infty} \left( \Psi_2 \right)$ by means of 
\begin{equation}
\langle O f_1  ,  f_2 \rangle_{\H_{1 / w_2}^{\infty} \left( \Psi_2 \right), \, \H_{w_2}^{1} \left( \Psi_2 \right)}
= \langle K  ,  f_1 \otimes f_2 \rangle_{\H_{1 / (w_1 \otimes w_2)}^{\infty,  \infty} \left( \Psi_1 \otimes \Psi_2 \right), \, \H_{w_1 \otimes w_2}^{1, 1} \left( \Psi_1 \otimes \Psi_2 \right)} 
\label{eq:Satz_3_1}
\end{equation}  for any $f_1 \in \H_{w_1}^{1} \left( \Psi_1 \right)$ and any $f_2 \in \H_{w_2}^{1} \left( \Psi_2 \right)$. Moreover,  
\begin{equation}
  \left\| K \right\|_{\H_{1 / (w_1 \otimes w_2)}^{\infty,   \infty} \left( \Psi_1 \otimes \Psi_2 \right)} \asymp \left\| O \right\|_{\H_{w_1}^{1} \left( \Psi_1 \right) \to \H_{1 / w_2}^{\infty} \left( \Psi_2 \right)} .
\label{eq:Satz_3_2}
\end{equation} 
Conversely, to any bounded linear  operator $O:\H_{w_1}^{1} \left( \Psi_1 \right)\to\H_{1 / w_2}^{\infty} \left( \Psi_2 \right)$ corresponds  a unique $K \in \H_{1 / (w_1 \otimes w_2)}^{\infty,   \infty} \left( \Psi_1 \otimes \Psi_2 \right)$,  called the \emph{kernel} of $O$, that satisfies \eqref{eq:Satz_3_1} for every $f_1 \in \H_{w_1}^{1} \left( \Psi_1 \right)$ and  $f_2 \in \H_{w_2}^{1} \left( \Psi_2 \right)$.
\end{theorem}

\proof\textbf{Step 1.} Note that $f_1 \otimes f_2 \in \H_{w_1 \otimes w_2}^{1,   1} \left( \Psi_1 \otimes \Psi_2 \right)$ as long as $f_1 \in \H_{w_1}^{1} \left( \Psi_1 \right)$ and $f_2 \in \H_{w_2}^{1} \left( \Psi_2 \right)$.  So, for $K \in \H_{1 / (w_1 \otimes w_2)}^{\infty, \, \infty} \left( \Psi_1 \otimes \Psi_2 \right)$,
\begin{align*}
  \big\vert \langle K   ,  f_1 \otimes f_2 \rangle &_{\H_{1 / (w_1 \otimes w_2)}^{\infty,  \infty} \left( \Psi_1 \otimes \Psi_2 \right),   \H_{w_1 \otimes w_2}^{1,   1} \left( \Psi_1 \otimes \Psi_2 \right)} \big\vert \notag \\
&   = \left\vert \sum_{i \in I} \sum_{j \in J} \big\langle K,  \widetilde \psi_{1,  \, i} \otimes \widetilde \psi_{2, \, j} \big\rangle  \left\langle \psi_{1,\,  i} \otimes \psi_{2,  \, j} , f_1 \otimes f_2 \right\rangle \right\vert \notag \\
&   \le \sum_{i \in I} \sum_{j \in J}  \big\vert \big\langle K,   \widetilde \psi_{1, \, i} \otimes \widetilde \psi_{2, \, j} \big\rangle  \left\langle \psi_{1, \, i} \otimes \psi_{2, \,  j}   ,   f_1\otimes f_2 \right\rangle  \big\vert  
\label{eq:Satz_3_1_Dualitaetserklaerung}
\\
&   \le\big\| C_{\widetilde{\Psi_1 \otimes \Psi_2}}  K \big\|_{\ell_{1 / (w_1 \otimes w_2)}^{\infty,  \infty}(I\times J)} \,  \big\| C_{\Psi_1 \otimes \Psi_2} (f_1 \otimes f_2 )\big\|_{\ell_{w_1 \otimes w_2}^{1,  1}(I\times J)}\notag
%\\
%&  = \left\| K \right\|_{\H_{1 / (w_1 \otimes w_2)}^{\infty,  \infty} \left( \Psi_1 \otimes \Psi_2 \right)} \, \big\| C_{\Psi_1} f_1\big\|_{\ell_{w_1 }^{  1}  } \big\|C_{\Psi_2}f_2 \big\|_{\ell_{ w_2}^{  1}  }
%\notag
\\
&  \asymp \left\| K \right\|_{\H_{1 / (w_1 \otimes w_2)}^{\infty, \, \infty} \left( \Psi_1 \otimes \Psi_2 \right)} \, \left\| f_1 \otimes f_2 \right\|_{\H_{w_1 \otimes w_2}^{1, \,  1} \left( \Psi_1 \otimes \Psi_2 \right)}
\notag
\\
&  < \infty,\notag
\end{align*}
where we used Proposition~\ref{Zerlegungssatz_Satz}~\ref{en:ii} for the second to last step.
%Therefore, the scalar product on the right hand side of (\ref{eq:Satz_3_1}) does indeed make sense.
Now, 
$$ 
\left( f_1 ,  f_2 \right) \mapsto \langle K   ,   f_1 \otimes f_2 \rangle_{\H_{1 / (w_1 \otimes w_2)}^{\infty,   \infty} \left( \Psi_1 \otimes \Psi_2 \right),   \H_{w_1 \otimes w_2}^{1, 1} \left( \Psi_1 \otimes \Psi_2 \right)}
$$
is a sesquilinear form on $\H_{w_1}^1(\Psi_1) \times \H_{w_2}^1(\Psi_2)$.
Therefore, for any fixed $f_1 \in \H_{w_1}^{1} \left( \Psi_1 \right)$, the  mapping $f_2 \mapsto \langle K   ,  f_1 \otimes f_2 \rangle_{\H_{1 / (w_1 \otimes w_2)}^{\infty, \, \infty} \left( \Psi_1 \otimes \Psi_2 \right),  \H_{w_1 \otimes w_2}^{1, \, 1} \left( \Psi_1 \otimes \Psi_2 \right)}$ is an antilinear functional on $\H_{w_2}^{1} \left( \Psi_2 \right)$ which we shall call $O f_1$. The mapping $f_1 \mapsto O f_1$ is linear and so (\ref{eq:Satz_3_1}) defines a linear operator $O:\H_{w_1}^{1} \left( \Psi_1 \right)\to \H_{1 / w_2}^{\infty} \left( \Psi_2 \right)$. This operator is bounded. Indeed,
\begin{align*}
 \big\vert \langle O f_1   ,  f_2 \rangle_{\H_{1 / w_2}^{\infty} \left( \Psi_2 \right),  \H_{w_2}^{1} \left( \Psi_2 \right)} \big\vert &
   = \left\vert \langle K   ,   f_1\otimes f_2 \rangle_{\H_{1 / (w_1 \otimes w_2)}^{\infty,   \infty} \left( \Psi_1 \otimes \Psi_2 \right),  \H_{w_1 \otimes w_2}^{1,   1} \left( \Psi_1 \otimes \Psi_2 \right)} \right\vert \notag
   \\
 & \leqslant \left\| K \right\|_{\H_{1 / (w_1 \otimes w_2)}^{\infty,  \infty} \left( \Psi_1 \otimes \Psi_2 \right)} 
      \left\| f_2 \right\|_{\H_{w_2}^{1} \left( \Psi_2 \right)} \left\| f_1 \right\|_{\H_{w_1}^{1} \left( \Psi_1 \right)},
\end{align*}
and therefore
\begin{equation*}
 \left\| O f_1 \right\|_{\H_{1 / w_2}^{\infty} \left( \Psi_2 \right)} 
  \leqslant \left\| K \right\|_{\H_{1 / (w_1 \otimes w_2)}^{\infty, \, \infty} \left( \Psi_1 \otimes \Psi_2 \right)} 
     \left\| f_1 \right\|_{\H_{w_1}^{1} \left( \Psi_1 \right)},
\label{eq:Satz_3_1_toto_2}
\end{equation*} 
and 
\begin{equation*}
 \left\| O \right\|_{\H_{w_1}^{1} \left( \Psi_1 \right) \to \H_{1 / w_2}^{\infty} \left( \Psi_2 \right)} 
  \leqslant \left\| K \right\|_{\H_{1 / (w_1 \otimes w_2)}^{\infty, \, \infty} \left( \Psi_1 \otimes \Psi_2 \right)}.
\label{eq:Satz_3_1_toto_3}
\end{equation*}
This proves one of the inequalities implied by (\ref{eq:Satz_3_2}) and indicates that the map $K \mapsto O$ is bounded. 

\textbf{Step 2.} Let us now assume that the bounded linear operator $O:\H_{w_1}^{1} \left( \Psi_1 \right)\to\H_{1 / w_2}^{\infty} \left( \Psi_2 \right)$ has two distinct   kernels $K_1$ and $K_2 \in \H_{1 / (w_1 \otimes w_2)}^{\infty,  \infty} \left( \Psi_1 \otimes \Psi_2 \right)$ satisfying \eqref{eq:Satz_3_1}. This implies in particular that for every $(i,j)\in I\times J$
%, let us assume that
%, for any $f_1 \in H_{w_1}^{1} \left( \Psi_1 \right)$ and any $f_2 \in H_{w_2}^{1} \left( \Psi_2 \right)$, and so for $\psi_{1, \, i} \in \Psi_1 \subset H_{w_1}^{1} \left( \Psi_1 \right)$ and $\psi_{2, \, j} \in \Psi_2 \subset H_{w_2}^{1} \left( \Psi_2 \right)$ %in particular,
\begin{equation*}
\big\langle O \widetilde\psi_{1, \, i} , \widetilde\psi_{2, \, j} \big\rangle_{\H_{1 / w_2}^{\infty} \left( \Psi_2 \right),   \H_{w_2}^{1} \left( \Psi_2 \right)}
= \big\langle K_1  ,   \widetilde\psi_{1, \, i} \otimes \widetilde\psi_{2, \, j} \big\rangle_{\H_{1 / (w_1 \otimes w_2)}^{\infty,  \infty} \left( \Psi_1 \otimes \Psi_2 \right),   \H_{w_1 \otimes w_2}^{1,  1} \left( \Psi_1 \otimes \Psi_2 \right)},
\label{eq:Satz_3_1_K_1}
\end{equation*} 
and
\begin{equation*}
\big\langle O \widetilde\psi_{1, \, i}  ,   \widetilde\psi_{2, \, j} \big\rangle_{\H_{1 / w_2}^{\infty} \left( \Psi_2 \right),  \H_{w_2}^{1} \left( \Psi_2 \right)}
= \big\langle K_2 , \widetilde\psi_{1,\,   i} \otimes \widetilde\psi_{2, \, j} \big\rangle_{\H_{1 / (w_1 \otimes w_2)}^{\infty,   \infty} \left( \Psi_1 \otimes \Psi_2 \right),  \H_{w_1 \otimes w_2}^{1,   1} \left( \Psi_1 \otimes \Psi_2 \right)}.
\label{eq:Satz_3_1_K_2}
\end{equation*} 
Hence, for every $(i,j)\in I\times J$ it holds
\begin{align*}
\big\langle K_1 ,   \widetilde\psi_{1, \, i} \otimes \widetilde\psi_{2, \, j} &\big\rangle_{\H_{1 / (w_1 \otimes w_2)}^{\infty,   \infty} \left( \Psi_1 \otimes \Psi_2 \right),   \H_{w_1 \otimes w_2}^{1,  1} \left( \Psi_1 \otimes \Psi_2 \right)} \notag
\\ &= \big\langle K_2   ,   \widetilde\psi_{1, \, i} \otimes \widetilde\psi_{2, \, j} \big\rangle_{\H_{1 / (w_1 \otimes w_2)}^{\infty,   \infty} \left( \Psi_1 \otimes \Psi_2 \right),   \H_{w_1 \otimes w_2}^{1,1} \left( \Psi_1 \otimes \Psi_2 \right)}.
\label{eq:Satz_3_1_K_f2_f1}
\end{align*}
From this assumption and the fact that, according to Proposition~\ref{Zerlegungssatz_Satz}~\ref{prop:(ii)}, any $F \in \H_{w_1 \otimes w_2}^{1, \, 1} \left( \Psi_1 \otimes \Psi_2 \right)$ can be expressed by an absolutely convergent series as $F=D_{\Psi_1\otimes\Psi_2}c$   we infer that
\begin{align*}
\langle K_1   ,   F  &\rangle_{\H_{1 / (w_1 \otimes w_2)}^{\infty, \, \infty} \left( \Psi_1 \otimes \Psi_2 \right),   \H_{w_1 \otimes w_2}^{1,  1} \left( \Psi_1 \otimes \Psi_2 \right)}
\\ &
    = \Big\langle K_1  , \sum_{i \in I} \sum_{j \in J} c_{i,j} \, \psi_{1,\,   i} \otimes \psi_{2,\,  j} \Big\rangle_{\H_{1 / (w_1 \otimes w_2)}^{\infty, \infty} \left( \Psi_1 \otimes \Psi_2 \right),  \H_{w_1 \otimes w_2}^{1,  1} \left( \Psi_1 \otimes \Psi_2 \right)} \\
 & = \sum_{i \in I} \sum_{j \in J} c_{i,  j} \,  \langle K_1   ,  \psi_{1, \,  i} \otimes \psi_{2, \, j} \rangle_{\H_{1 / (w_1 \otimes w_2)}^{\infty,   \infty} \left( \Psi_1 \otimes \Psi_2 \right),   \H_{w_1 \otimes w_2}^{1,   1} \left( \Psi_1 \otimes \Psi_2 \right)} \\
 & = \sum_{i \in I} \sum_{j \in J} c_{i,   j}  \langle K_2   ,   \psi_{1,\,  i} \otimes \psi_{2,\,   j} \rangle_{\H_{1 / (w_1 \otimes w_2)}^{\infty, \, \infty} \left( \Psi_1 \otimes \Psi_2 \right),  \H_{w_1 \otimes w_2}^{1,  1} \left( \Psi_1 \otimes \Psi_2 \right)} \\
 &  = \Big\langle K_2   ,   \sum_{i \in I} \sum_{j \in J} c_{i, j} \, \psi_{1,  \, i} \otimes \psi_{2, \, j} \Big\rangle_{\H_{1 / (w_1 \otimes w_2)}^{\infty,   \infty} \left( \Psi_1 \otimes \Psi_2 \right),   \H_{w_1 \otimes w_2}^{1,  1} \left( \Psi_1 \otimes \Psi_2 \right)} \\
 & = \langle K_2   ,   F \rangle_{\H_{1 / (w_1 \otimes w_2)}^{\infty,  \infty} \left( \Psi_1 \otimes \Psi_2 \right), \H_{w_1 \otimes w_2}^{1,  1} \left( \Psi_1 \otimes \Psi_2 \right)},
%\label{eq:Satz_3_1_K_f}
\end{align*}
and so $K_1 = K_2$ according to Proposition~\ref{Zerlegungssatz_Satz}~\ref{en:dual}.  In other words, the map $K \mapsto O$ is injective. 

\textbf{Step 3.} Let us now assume that the linear operator $O:\H_{w_1}^{1} \left( \Psi_1 \right)\to \H_{1 / w_2}^{\infty} \left( \Psi_2 \right)$ is bounded and consider %what is often referred to as 
its Galerkin matrix, that is  the matrix $k$ whose $(i,   j)$-th element is defined by
\begin{equation}\label{eq:k}
k_{i ,   j} := \big\langle O \widetilde \psi_{1,  i}  , \widetilde \psi_{2,   j} \big\rangle_{\H_{1/w_2}^{\infty} \left( \Psi_2 \right),   \H_{w_2}^{1} \left( \Psi_2 \right)},\quad (i,   j) \in I \times J.
\end{equation} 
From \eqref{eq:k}, and \eqref{eq:bound-norm-frame-elements} we infer that
\begin{equation*}
\begin{split}
 \left\vert k_{i , \, j} \right\vert &
   \leqslant \big\| O \big\|_{\H_{w_1}^{1} \left( \Psi_1 \right) \to \H_{1 / w_2}^{\infty} \left( \Psi_2 \right)} \, \big\| \widetilde \psi_{1, \, i} \big\|_{\H_{w_1}^{1} \left( \Psi_1 \right)} \,
        \big\| \widetilde \psi_{2, \, j} \big\|_{\H_{w_2}^{1} \left( \Psi_2 \right)} \\
  & \lesssim \left\| O \right\|_{\H_{w_1}^{1} \left( \Psi_1 \right) \to \H_{1 / w_2}^{\infty} \left( \Psi_2 \right)} \, \omega_{1, \, i} \, \omega_{2, \, j}.
 \end{split}
\label{eq:Norm_k_O}
\end{equation*} 
In other words, $k \in \ell_{1/(w_1 \otimes w_2)}^{\infty, \, \infty} (I \times J)$ as long as the operator $O:\H_{w_1}^{1} \left( \Psi_1 \right)\to\H_{1 / w_2}^{\infty} \left( \Psi_2 \right)$ is bounded. Subsequently, we show that   $k$ satisfies \eqref{eq:Correspondence_2}, which, according to Lemma~\ref{thm3}, implies that there is a unique $K \in \H_{1/(w_1 \otimes w_2)}^{\infty, \, \infty} \left( \Psi_1 \otimes \Psi_2 \right)$ that
%, together with $k \in \ell_{1/(w_1 \otimes w_2)}^{\infty, \, \infty} (I \times J)$, 
satisfies $C_{\widetilde{\Psi_1\otimes \Psi_2}}K=k$. Indeed, for any $(i, j) \in I \times J$,
\begin{align*}
\Big( C_{\widetilde{\Psi_1 \otimes \Psi_2}} & D_{\Psi_1 \otimes \Psi_2} k \, \Big)_{r, \, s} \notag
\\
& = \sum_{i \in I} \sum_{j \in J} \big\langle O \widetilde\psi_{1, \, i}   ,   \widetilde\psi_{2, \, j} \big\rangle_{\H_{1/w_2}^{\infty} \left( \Psi_2 \right),  \H_{w_2}^{1} \left( \Psi_2 \right)} \, \big\langle \psi_{1, \, i} \otimes \psi_{2, \, j}   ,  \widetilde\psi_{1, \, r} \otimes \widetilde\psi_{2, \, s} \big\rangle_{\H_1\otimes\H_2} \notag
\\
 & = \sum_{i \in I} \sum_{j \in J} \big\langle O \widetilde\psi_{1, \, i}    ,   \widetilde\psi_{2, \, j} \big\rangle_{\H_{1/w_2}^{\infty} \left( \Psi_2 \right),\H_{w_2}^{1} \left( \Psi_2 \right)} \, \overline{\big\langle \psi_{1, \, i}    ,  \widetilde\psi_{1, \, r}  \big\rangle_{\H_1 } }
 \, \big\langle \psi_{2, \, j}  ,  \widetilde\psi_{2, \, s} \big\rangle_{ \H_2}  \notag
 \\
 & = \sum_{i \in I} \langle O \widetilde\psi_{1,  \, i}    ,  \widetilde\psi_{2, \, s} \rangle_{\H_{1/w_2}^{\infty} \left( \Psi_2 \right),   \H_{w_2}^{1} \left( \Psi_2 \right)} \,
\overline{\big\langle \psi_{1, \, i}    ,  \widetilde\psi_{1, \, r}  \big\rangle_{ \H_1} }
 \\
 & = \sum_{i \in I} \langle \widetilde\psi_{1, \, i}   ,   O' \widetilde\psi_{2, \, s} \rangle_{\H_{w_1}^{1} \left( \Psi_1 \right),  \H_{1/w_1}^{\infty} \left( \Psi_1 \right)} \,
 \big\langle   \widetilde\psi_{1, \, r} ,\psi_{1, \, i}     \big\rangle_{ \H_1}   \notag
 \\
 & = \langle \widetilde\psi_{1, \, r}   ,   O' \widetilde\psi_{2, \, s} \rangle_{\H_{w_1}^{1} \left( \Psi_1 \right), \, \H_{1/w_1}^{\infty} \left( \Psi_1 \right)} \notag \\
 & = \langle O \widetilde\psi_{1, \, r}    ,   \widetilde\psi_{2, \, s} \rangle_{\H_{1/w_2}^{\infty} \left( \Psi_2 \right),  \H_{w_2}^{1} \left( \Psi_2 \right)} = k_{r , \, s}, \notag
\label{eq:Satz_3_1_toto_3}
\end{align*}
where $O'$ stands for the operator that satisfies
\begin{equation*}
   \langle O f_1    ,  f_2  \rangle_{\H_{1/w_2}^{\infty} \left( \Psi_2 \right),   \H_{w_2}^{1} \left( \Psi_2 \right)} = \langle f_1   ,   O' f_2 \rangle_{\H_{w_1}^{1} \left( \Psi_1 \right),   \H_{1/w_1}^{\infty} \left( \Psi_1 \right)},
\label{eq:paeadjungierter_Operator}
\end{equation*}
for every $f_1 \in \H_{w_1}^{1} \left( \Psi_1 \right)$ and  $f_2 \in \H_{w_2}^{1} \left( \Psi_2 \right)$. The change of the order of summation over $i$ and $j$ that we made here is justified by the fact that
\begin{equation*}
\left\{ \big\langle O \widetilde\psi_{1, \, i}    ,  \widetilde\psi_{2, \, j} \big\rangle_{\H_{1/w_2}^{\infty} \left( \Psi_2 \right),   \H_{w_2}^{1} \left( \Psi_2 \right)} \,   \overline{\big\langle \psi_{1, \, i}    ,  \widetilde\psi_{1, \, r}  \big\rangle_{ \H_1} } \, \big\langle \psi_{2, \, j}   ,  \widetilde\psi_{2, \, s} \big\rangle_{ \H_2}  \right\}_{(i, \, j) \in I^2} \in \ell^{1,  1} (I \times J)
\label{eq:toto}
\end{equation*}
for any given $(r, \, s) \in I \times J$, as
\begin{equation*}
\left\{ \big\langle O \widetilde\psi_{1, \, i}    ,   \widetilde\psi_{2, \, j} \big\rangle_{\H_{1/w_2}^{\infty} \left( \Psi_2 \right),   \H_{w_2}^{1} \left( \Psi_2 \right)} \right\}_{(i, \, j) \in I \times J} = \{k_{i,j}\}_{(i,j)\in I\times J} \in \ell_{1/(w_1 \otimes w_2)}^{\infty,  \infty}  (I \times J)
\label{eq:toto_2}
\end{equation*}
and
\begin{align*}
\left\{\overline{\big\langle \psi_{1, \, i}    ,  \widetilde\psi_{1, \, r}  \big\rangle_{ \H_1} }  \big\langle \psi_{2, \, j}   ,   \widetilde\psi_{2, \, s} \big\rangle_{ \H_2}  \right\}_{(i, \, j) \in I \times J}
=\left\{\overline{ (C_{ \Psi_1 }\widetilde\psi_{1, \, r})_i}(C_{ \Psi_2}    \widetilde \psi_{2, \, s})_j\right\}& _{(i,j)\in I\times J}
\\& \hspace{-8pt}\in \ell_{w_1 \otimes w_2}^{1,  1} (I \times J),
\label{eq:toto_3}
\end{align*}
by Proposition~\ref{prop:properties}~\ref{bp:6} and  Lemma~\ref{lem:bound-norm-frame-elements}.
 According to what we already proved, to this unique $K \in \H_{1 / (w_1 \otimes w_2)}^{\infty,  \infty} \left( \Psi_1 \otimes \Psi_2 \right)$ corresponds a linear operator mapping $\H_{w_1}^{1} \left( \Psi_1 \right)$ to $\H_{1 / w_2}^{\infty} \left( \Psi_2 \right)$ via \eqref{eq:Satz_3_1}. This operator is nothing, but the bounded linear operator $O$. Indeed, let us assume that there is more than one linear operator mapping $\H_{w_1}^{1} \left( \Psi_1 \right)$ to  $\H_{1 / w_2}^{\infty} \left( \Psi_2 \right)$, say,  $O_1$ and  $O_2$, that correspond to the kernel $K \in \H_{1 / (w_1 \otimes w_2)}^{\infty,  \infty} \left( \Psi_1 \otimes \Psi_2 \right)$  via \eqref{eq:Satz_3_1}. 
% In other words, let us assume that any element of $ \H_{w_1 \otimes w_2}^{1,   1} \left( \Psi_1 \otimes \Psi_2 \right)$ satisfies  \eqref{eq:Satz_3_1}. 
So in particular for  $\psi_{1, \, i} \otimes \psi_{2, \, j}\in  \H_{w_1 \otimes w_2}^{1,   1} \left( \Psi_1 \otimes \Psi_2 \right)$ we get %where $\psi_{1,  \, i} \in \Psi_1 \subset H_{w_1}^{1} \left( \Psi_1 \right)$ and $\psi_{2, \, j} \in \Psi_2 \subset H_{w_2}^{1} \left( \Psi_2 \right)$ in particular satisfies \eqref{eq:Satz_3_1}, i.e.
\begin{equation*}
\big\langle O_1 \widetilde\psi_{1, \, i}   ,   \widetilde\psi_{2, \, j} \big\rangle_{\H_{1 / w_2}^{\infty} \left( \Psi_2 \right),  \H_{w_2}^{1} \left( \Psi_2 \right)}
= \big\langle K   ,   \widetilde\psi_{1, \, i} \otimes \widetilde\psi_{2, \, j} \big\rangle_{\H_{1 / (w_1 \otimes w_2)}^{\infty,  \infty} \left( \Psi_1 \otimes \Psi_2 \right),   \H_{w_1 \otimes w_2}^{1,   1} \left( \Psi_1 \otimes \Psi_2 \right)},
\label{eq:Satz_3_1_O_1}
\end{equation*} 
and
\begin{equation*}
\big\langle O_2 \widetilde\psi_{1, \, i}  ,  \widetilde\psi_{2, \, j} \big\rangle_{\H_{1 / w_2}^{\infty} \left( \Psi_2 \right),  \H_{w_2}^{1} \left( \Psi_2 \right)}
= \big\langle K  ,   \widetilde\psi_{1, \, i} \otimes \widetilde\psi_{2, \, j} \big\rangle_{\H_{1 / (w_1 \otimes w_2)}^{\infty,   \infty} \left( \Psi_1 \otimes \Psi_2 \right),  \H_{w_1 \otimes w_2}^{1,  1} \left( \Psi_1 \otimes \Psi_2 \right)}.
\label{eq:Satz_3_1_O_2}
\end{equation*} 
Thus
\begin{align*}
\big\langle O_1 \widetilde\psi_{1, \, i}   ,   \widetilde\psi_{2, \, j} &\big\rangle_{\H_{1 / w_2}^{\infty} \left( \Psi_2 \right),   \H_{w_2}^{1} \left( \Psi_2 \right)}
= \big\langle O_2 \widetilde\psi_{1, \, i}   ,   \widetilde\psi_{2, \, j} \big\rangle_{\H_{1 / w_2}^{\infty} \left( \Psi_2 \right),   \H_{w_2}^{1} \left( \Psi_2 \right)}.
%\label{eq:Satz_3_1_O_f2_f1}
\end{align*}
This equality, and the fact that, according to Proposition~\ref{thm0}, any $f_1 \in \H_{w_1}^{1} \left( \Psi_1 \right)$ and any $f_2 \in \H_{w_2}^{1} \left( \Psi_2 \right)$ can be expressed by an absolutely convergent series, we infer that 
\begin{align*}
\langle O_1 f_1  ,   f_2 \rangle_{\H_{1 / w_2}^{\infty} \left( \Psi_2 \right),  \H_{w_2}^{1} \left( \Psi_2 \right)} &
  = \Big\langle O_1 \sum_{i \in I} c_{1, \, i} \, \widetilde\psi_{1, \, i}  ,   \sum_{j \in J} c_{2, \, j} \, \widetilde\psi_{2, \, j} \Big\rangle_{\H_{1 / w_2}^{\infty} \left( \Psi_2 \right),   \H_{w_2}^{1} \left( \Psi_2 \right)} \notag 
  \\
 & = \sum_{i \in I} \sum_{j \in J} c_{1, \, i}  \, c_{2, \, j}   \big\langle O_1  \widetilde\psi_{1, \, i}   ,   \widetilde\psi_{2, \, j} \big\rangle_{\H_{1 / w_2}^{\infty} \left( \Psi_2 \right),   \H_{w_2}^{1} \left( \Psi_2 \right)}\notag 
 \\
 & = \sum_{i \in I} \sum_{j \in J} c_{1, \, i}  \, c_{2, \, j} \, \big\langle O_2   \widetilde\psi_{1, \, i}  ,  \widetilde\psi_{2, \, j} \big\rangle_{\H_{1 / w_2}^{\infty} \left( \Psi_2 \right),   \H_{w_2}^{1} \left( \Psi_2 \right)} 
 \\
 & = \Big\langle O_2 \sum_{i \in I} c_{1, \, i} \, \widetilde\psi_{1, \, i}   ,   \sum_{j \in J} c_{2, \, j} \, \widetilde\psi_{2, \, j} \Big\rangle_{\H_{1 / w_2}^{\infty} \left( \Psi_2 \right),   \H_{w_2}^{1} \left( \Psi_2 \right)} \notag
 \\
 & = \langle O_2 f_1   , f_2 \rangle_{\H_{1 / w_2}^{\infty} \left( \Psi_2 \right),   \H_{w_2}^{1} \left( \Psi_2 \right)}.\notag
\label{eq:Satz_3_1_toto1}
\end{align*} 
Thus $O_1 f_1 = O_2 f_1$ for every $f_1 \in \H_{w_1}^{1} \left( \Psi_1 \right)$  or, in other words, $O_1 = O_2$. Consequently, the mapping $O \mapsto K$ that takes the space $B \big( \H_{w_1}^{1}  ( \Psi_1  ),  \H_{1 / w_2}^{\infty} \left( \Psi_2 \right)\big)$ to  $\H_{1 / (w_1 \otimes w_2)}^{\infty,  \infty} \left( \Psi_1 \otimes \Psi_2 \right)$, is surjective and invertible and, as we proved before, bounded and injective. Therefore, its inverse is also bounded according to the inverse mapping theorem, which is expressed by the first inequality in \eqref{eq:Satz_3_2}. \pbox

The so-called \emph{tensor product property}, discovered by Feichtinger \cite[Theorem 7D]{Feichtinger_1981} (see also \cite{Losert}), was generalised for co-orbit spaces associated with integrable group representations in \cite{Balazs_2019}. In what follows, we shall obtain a similar result for co-orbit spaces associated with $\mathcal{A}$-localised frames.
Let us recollect the definition of the \emph{projective tensor product} of two Banach spaces $B_1$ and $B_2$,
given by
\begin{equation*}
B_1 \widehat{\otimes}_\pi B_2 
   := \left\lbrace F  =  \sum_{r \in R}f_{ r}  \otimes  g_{  r}:\quad     f_{  r} \in B_1, \, g_{  r}\in B_2,  \quad \sum_{r \in R}  \| f_{ r} \|_{B_1} \,  \| g_{   r} \|_{B_2} < \infty  \right\rbrace.
\label{eq:projektives_Tensorprodukt_Ausdruck}
\end{equation*} 
The function
\begin{equation}
\begin{split}
\| F \|_{B_1 \widehat{\otimes}_\pi B_2} := \inf \left( \sum_{r \in R}  \| f_{ r} \|_{B_1} \,  \| g_{r} \|_{B_2} \right),
\end{split}
\label{eq:projektives_Tensorproduktraumnorm}
\end{equation} 
the infimum being taken over all representations
$\sum_{r \in R}   f_{  r}  \otimes   g_{r}$ of $F \in B_1 \widehat{\otimes}_\pi B_2$,   defines a norm for $B_1 \widehat{\otimes}_\pi B_2$.

\begin{theorem}\label{Isomorphism} Let   $\Psi_n\hspace{-1pt}\subset\hspace{-1pt} \H_n$ be an $\mathcal{A}_n$-localised frame, and $w_n$ be an $\mathcal{A}_n$-admissible weight for $n=1,2$. Then
\begin{equation}
\begin{split}
\H_{w_1 \otimes w_2}^{1,  1} \left( \Psi_1 \otimes \Psi_2 \right) =
  \H_{w_1}^1 \left( \Psi_1 \right) \widehat{\otimes}_\pi \H_{w_2}^1 \left( \Psi_2 \right)
\end{split}
\label{eq:Isomorphism}
\end{equation}
with equivalent norms. 
\end{theorem}

\proof %First of all, we note that using Proposition~\ref{prop:properties}~(vi) and arguments similar to those that we used to obtain \eqref{eq:Norm_psi_1} and \eqref{eq:Norm_psi_2}   allow us to conclude that for every $i\in I$ and $j\in J$
%\begin{equation}\label{eq:Norm_psi_1_bis}
% \| \psi_{1, \, i} \|_{\H_{w_1}^1 \left( \Psi_1 \right)} \lesssim w_{1, \, i},\quad \text{and}\quad
 %\| \psi_{2, \, j} \|_{\H_{w_2}^1 \left( \Psi_2 \right)} \lesssim w_{2, \, j},
%\end{equation} for any $j \in J$.
 Let  $F \in \H_{w_1 \otimes w_2}^{1,1} \left( \Psi_1 \otimes \Psi_2 \right)$. Then, according to Proposition~\ref{Zerlegungssatz_Satz}~\ref{prop:(ii)}, the function $F$ can be expressed as $F=D_{\Psi_1\otimes \Psi_2}c$ for some $c\in\ell^{1,1}_{w_1\otimes w_2}(I\times J)$ such that $\|c\|_{\ell^{1,1}_{w_1\otimes w_2}(I\times J)}\asymp \|F\|_{\H^{1,1}_{w_1\otimes w_2}(\Psi_1\otimes \Psi_2)}$. Combining this with \eqref{eq:bound-norm-frame-elements}   results in 
\begin{equation*}
\begin{split}
\| F \|_{\H_{w_1}^1 \left( \Psi_1 \right)  \widehat{\otimes}_\pi  \H_{w_2}^1 \left( \Psi_2 \right)} &
  \leqslant \sum_{i \in I} \sum_{j \in J} | c_{i,  j } |\,\|  \psi_{1, \, i} \|_{\H_{w_1}^1 \left( \Psi_1 \right)}    \| \psi_{2, \, j} \|_{\H_{w_2}^1 \left( \Psi_2 \right)} \\
  &   \lesssim \sum_{i \in I} \sum_{j \in J} \vert c_{i,  j} \vert \, w_{1, \, i} \, w_{2, \, j} \\
  & =  \left\| c \right\|_{\ell_{w_1 \otimes w_2}^{1,   1} \left( I \times J \right)} \asymp \| F \|_{\H_{w_1 \otimes w_2}^{1,   1} \left( \Psi_1 \otimes \Psi_2 \right)} .
\end{split}
\label{eq:Isomorphism_toto}
\end{equation*} 
Thus $\H_{w_1 \otimes w_2}^{1,   1} \left( \Psi_1 \otimes \Psi_2 \right)$ is continuously embedded in $\H_{w_1}^1 \left( \Psi_1 \right)   \widehat{\otimes}_\pi  \, \H_{w_2}^1 \left( \Psi_2 \right)$.

 Let now $F \in \H_{w_1}^1 \left( \Psi_1 \right)   \widehat{\otimes}_\pi\,   \H_{w_2}^1 \left( \Psi_2 \right)$. In other words, let
$F  =  \sum_{r \in R}  f_{r}   \otimes  g_{r}$ 
 where $ f_{ r} \in \H_{w_1}^1 \left( \Psi_1 \right)$ and  $ g_{r} \in \H_{w_2}^1 \left( \Psi_2 \right)$ satisfy 
$$ 
\sum_{r \in R}  \| f_{ r} \|_{\H_{w_1}^1 \left( \Psi_1 \right)} \,  \| g_{r} \|_{\H_{w_2}^1 \left( \Psi_2 \right)} < \infty.
$$
Then
\begin{equation*}
\begin{split}
\| F \|_{\H_{w_1 \otimes w_2}^{1,  1} \left( \Psi_1 \otimes \Psi_2 \right)} &
  = \| C_{\widetilde{\Psi_1 \otimes \Psi_2}}F\|_{\ell_{w_1 \otimes w_2}^{1,  1} \left( I \times J \right)}  \\
  & = \sum_{i \in I} \sum_{j \in J} \Big\vert   \sum_{r \in R} \big(C_{\widetilde{\Psi_1 \otimes \Psi_2}} ( f_{  r}  \otimes  g_{r} )\big)_{i,j}   \Big\vert  \, w_{1, \, i} \,  w_{2, \,j}\\ 
  & \leqslant \sum_{r \in R}   \Big( \sum_{i \in I} \big\vert  \big(C_{\widetilde{\Psi}_1}f_{ r} \big)_i  \big\vert  w_{1, \, i} \Big) \Big( \sum_{j \in J} \big\vert \big(C_{\widetilde{\Psi}_2} g_{r}  \big)_j \big\vert   w_{2, \,j} \Big)   \\
  & = \sum_{r \in R}   \| C_{\widetilde{\Psi}_1}   f_{r} \|_{\ell_{w_1}^1 \left( I \right)} \, \| C_{\widetilde{\Psi}_2}  g_{r} \|_{\ell_{w_2}^1 \left( J \right)} \\
  & = \sum_{r \in R}  \| f_{  r} \|_{\H_{w_1}^1 \left( \Psi_1 \right)} \, 
    \| g_{r} \|_{\H_{w_2}^1 \left( \Psi_2 \right)}  < \infty  .
\end{split}
\label{eq:Isomorphism_toto_2}
\end{equation*}

The change of the order of summation over $i,   j, $ and $  r$  above is justified by the absolute convergence of the series. Taking the infimum over all possible representations of $F$ in the inequality above then shows that  $\H_{w_1}^1 \left( \Psi_1 \right) \widehat{\otimes}_\pi \, \H_{w_2}^1 \left( \Psi_2 \right)$ is also continuously embedded in $\H_{w_1 \otimes w_2}^{1, 1} \left( \Psi_1 \otimes \Psi_2 \right)$, which concludes the proof. \pbox

The previous theorem allows us to characterise  the class of operators with kernels in $\H^{1,1}_{w_1\otimes w_2}(\Psi_1\otimes\Psi_2)$.
\begin{comment}
Since  along the chain of coorbit spaces this corresponds to the innermost class of operator kernels, this result is often referred to in the literature as the inner kernel theorem \cite{Feichtinger_Jakobsen}.
\end{comment}

 \begin{theorem}[Inner Kernel Theorem]\label{innere_kernsatz} Let $\Psi_n$ be an $\mathcal{A}_n$-localised frame for a Hilbert space $\H_1$,  $w_n$ an $\mathcal{A}_n$-admissible weight $n=1,2$. 
 Then there is an isomorphism between $\H^{1,1}_{w_1\otimes w_2}(\Psi_1\otimes \Psi_2)$ and the space
\begin{align*}
\mathcal{B}:=  \Big\lbrace &O\in  B\big(\H^\infty_{1/w_1}(\Psi_1),\H^1_{w_2}(\Psi_2)\big):\ \text{there exist } f_{  r} \in H^1_{w_1}(\Psi_1),g_r\in H^1_{w_2}(\Psi_2)  \text{ s.t.}\\ & O  =  \sum_{r \in R}\langle\, \cdot\,  ,f_{ r}\rangle_{H^\infty_{1/w_1}(\Psi_1),H^1_{w_1}(\Psi_1)} \,  g_{  r},   \text{ and }   \, \sum_{r \in R}  \| f_{ r} \|_{H^1_{w_1}(\Psi_1)} \,  \| g_{   r} \|_{H^1_{w_2}(\Psi_2)} < \infty  \Big\rbrace.
\end{align*}
In particular, if any of the two,  
$$
K\in \H^{1,1}_{w_1\otimes w_2}( \Psi_1\otimes \Psi_2) \ \text{ or }\ O\in\mathcal{B},$$
is given, the other is uniquely determined by means of
$$
\langle K,f_1\otimes f_2\rangle_{\H^{1,1}_{w_1\otimes w_2}(\Psi_1\otimes\Psi_2),\H^{\infty,\infty}_{1/(w_1\otimes w_2)}(\Psi_1\otimes\Psi_2)}=\langle Of_1,f_2\rangle_{\H^{1}_{ w_2}(\Psi_2),\H^{\infty}_{1/ w_2}(\Psi_2)}
$$
 for any $f_1\in \H^{\infty}_{1/ w_1}(\Psi_1)$ and any $f_2\in \H^{\infty}_{1/ w_2}(\Psi_2)$.
 \end{theorem}
 \proof First of all, we note that the condition
 \begin{equation}
     \sum_{r \in R}  \| f_{ r} \|_{H^1_{w_1}(\Psi_1)} \,  \| g_{   r} \|_{H^1_{w_2}(\Psi_2)} < \infty
\label{eq:Endlichkeit}
 \end{equation} implies that $O\in  B\big(\H^\infty_{1/w_1}(\Psi_1),\H^1_{w_2}(\Psi_2)\big)$. Now, to prove the result, it suffices to establish an isomorphism between $\mathcal{B}$ and $\H^1_{w_1}(\Psi_1)\widehat{\otimes}_\pi \H^1_{w_2}(\Psi_2)$ and apply Theorem~\ref{Isomorphism}, according to which $\H^1_{w_1}(\Psi_1)\widehat{\otimes}_\pi \H^1_{w_2}(\Psi_2)$ and $\H_{w_1 \otimes w_2}^{1,  1} \left( \Psi_1 \otimes \Psi_2 \right)$ are isomorphic.
 \begin{comment}
 It thus suffices to establish an isomorphism between 
 \begin{align*}
     \mathcal{A}':=  \Big\lbrace   O  =  \sum_{r \in R}\langle\, \cdot\,  ,f_{ r}\rangle_{H^\infty_{1/w_1}(\Psi_1),H^1_{w_1}(\Psi_1)} &\,  g_{  r}:\  f_{  r} ,g_r\in H^1_{w_1}(\Psi_1), \\ &\text{ and }   \, \sum_{r \in R}  \| f_{ r} \|_{H^1_{w_1}(\Psi_1)} \,  \| g_{   r} \|_{H^1_{w_1}(\Psi_1)} < \infty  \Big\rbrace
     \end{align*}
 and $\H^1_{w_1}(\Psi_1)\widehat{\otimes}\H^1_{w_2}(\Psi_2)$, and apply Theorem~\ref{Isomorphism}.
  \end{comment}
Let us define the mapping 
$\mathcal{I}:\H^1_{w_1}(\Psi_1)\widehat{\otimes}_\pi \H^1_{w_2}(\Psi_2)\to \mathcal{B}$ by 
$$
\mathcal{I}\Big(\sum_{r\in R}f_r\otimes g_r\Big)=\sum_{r\in R}\langle\,\cdot\, ,f_r\rangle_{\H^\infty_{1/w_1}(\Psi_1),\H^1_{w_1}(\Psi_1)}\, g_r,
$$
and first show that it is well-defined.  
For any two representations $\sum_{r\in R}f_r\otimes g_r$  and $\sum_{s\in S}h_s\otimes k_s$ of $H \in \H^1_{w_1}(\Psi_1)\widehat\otimes_\pi \H^1_{w_2}(\Psi_2)$, any   $f_1\in\H^\infty_{1/w_1}(\Psi_1)$ and any $f_2\in\H^\infty_{1/w_2}(\Psi_2)$,
\begin{align*}
\Big\langle \sum_{r\in R}f_r\otimes g_r,f_1\otimes f_2&\Big\rangle_{\H^{1,1}_{w_1\otimes w_2}(\Psi_1\otimes\Psi_2),\H^{\infty,\infty}_{1/(w_1\otimes w_2)}(\Psi_1\otimes\Psi_2)} \\&=\Big\langle \sum_{s\in S}h_s\otimes k_s,f_1\otimes f_2\Big\rangle_{\H^{1,1}_{w_1\otimes w_2}(\Psi_1\otimes\Psi_2),\H^{\infty,\infty}_{1/(w_1\otimes w_2)}(\Psi_1\otimes\Psi_2)}.
\end{align*}
Changing the order of scalar products and summation results in
 \begin{align}\label{eq:I(H1)=I(H2)}
  \sum_{r\in R}\langle f_1,f_r\rangle_{\H^\infty_{1/w_1}(\Psi_1),\H^1_{w_1}(\Psi_1)}&\langle g_r,f_2\rangle_{\H^1_{w_2}(\Psi_2),\H^\infty_{1/w_2}(\Psi_2)}  \\& =  \sum_{s\in S}\langle f_1,h_s\rangle_{\H^\infty_{1/w_1}(\Psi_1),\H^1_{w_1}(\Psi_1)}\langle k_s,f_2\rangle_{\H^1_{w_2}(\Psi_2),\H^\infty_{1/w_2}(\Psi_2)}  \notag
 \end{align} for any   $f_1\in\H^\infty_{1/w_1}(\Psi_1)$ and any $f_2\in\H^\infty_{1/w_2}(\Psi_2)$,
 from which we infer that $\mathcal{I}\big(\sum_{r\in R}f_r\otimes g_r\big)=\mathcal{I}\big(\sum_{s\in S}h_s\otimes k_s\big)$. %Hahn-Banach sepataration theorem

It is easy to see that  $\mathcal
 I$ is surjective and it therefore remains to prove that is also injective. Let the elements $H_1=\sum_{r\in R}f_r\otimes g_r$ and $  H_2=\sum_{s\in S}h_s\otimes k_s$ of $\H^1_{w_1}(\Psi_1)\widehat\otimes_\pi \H^1_{w_2}(\Psi_2)$ be such that $\mathcal{I}(H_1)=\mathcal{I}(H_2)$. In other words, \eqref{eq:I(H1)=I(H2)} holds for any $f_1\in\H^\infty_{1/w_1}(\Psi_1)$ and any $f_2\in\H^\infty_{1/w_2}(\Psi_2)$.
 From this and \eqref{eq:Endlichkeit}, which holds for any element of $\H^1_{w_1}(\Psi_1)\widehat\otimes_\pi \H^1_{w_2}(\Psi_2)$, we deduce that
 \begin{align*}
 \langle H_1,f_1\otimes f_2\rangle&_{\H^{1,1}_{w_1\otimes w_2}(\Psi_1\otimes \Psi_2),\H^{\infty,\infty}_{1/(w_1\otimes w_2)}(\Psi_1\otimes \Psi_2)}
 \\
 &\hspace{2cm}= \langle H_2,f_1\otimes f_2\rangle_{\H^{1,1}_{w_1\otimes w_2}(\Psi_1\otimes \Psi_2),\H^{\infty,\infty}_{1/(w_1\otimes w_2)}(\Psi_1\otimes \Psi_2)}.
 \end{align*}
Since  $D_{\Psi_1\otimes\Psi_2}:\ell^{\infty,\infty}_{1/(w_1\otimes w_2}(I\times J)\to \H^{\infty,\infty}_{1/(w_1\otimes w_2}(\Psi_1\otimes \Psi_2)$ is surjective (Lemma~\ref{thm3}) 
and converges weak-$\ast$ unconditionally  (Proposition~\ref{Zerlegungssatz_Satz}~\ref{prop:(i)}) we may consider linear combinations of $f_1\otimes f_2=\psi_{1,\, i}\otimes \psi_{2,\, j}$ to show that  $H_1=H_2$ in $\H^{1,1}_{w_1\otimes w_2}(\Psi_1\otimes\Psi_2).$ According to Theorem~\ref{Isomorphism}, this equality also holds in $\H^1_{w_1}(\Psi_1)\widehat\otimes_\pi \H^1_{w_2}(\Psi_2)$, which concludes the proof.
 \pbox

Now that we established the correspondence between operators and kernels in the extremes of the chain of co-orbit spaces, we investigate whether it is possible to characterise certain intermediate classes of operators via the properties of their kernels, in particular in terms of certain mixed-norm co-orbit spaces. The arguments crucially rely on a version of Schur's test (see, for example, \cite[Propositions 5.2 and 5.4]{tao} or \cite{maddox}) and the representation of an operator by its Galerkin  matrix defined in \eqref{eq:k} by  $$k_{i,j}:=\big\langle O\widetilde{\psi}_{1,\, i},\widetilde{\psi}_{2,\, j}\big\rangle_{\H^\infty_{1/w_2}(\Psi_2),\H_{w_2}^1(\Psi_2)}.$$ We note that the characterisation given by this matrix was established in \cite[Proposition~6]{Balazs_2017}. What is new here is that we relate the Galerkin matrix representing an operator to its kernel.

\begin{theorem}\label{Supergeiler_Satz} Let $1\le p,q \le \infty$,  $\frac{1}{p} + \frac{1}{q} =1$, $\Psi_n\subset \H_n$ be an $\mathcal{A}_n$-localised frame, and $w_n$ be an $\mathcal{A}_n$-admissible weight for $n=1,2$. %Moreover, let $O$ be a bounded linear  operator mapping $\H_{w_1}^1 (\Psi_1)$ to $\H_{1/w_2}^{\infty} (\Psi_2)$   with kernel $K \in \H_{1 / (w_1 \otimes w_2)}^{\infty, \, \infty} \left( \Psi_1 \otimes \Psi_2 \right)$.
Then
\begin{enumerate}[label=(\roman*)]
    \item\label{thm:super-i} A linear operator $O:\H_{w_1}^1 (\Psi_1)\to\H_{1/w_2}^p (\Psi_2)$ is  bounded  if and only if its kernel $K$ is contained in $\mathfrak{H}_{1/(w_1 \otimes w_2)}^{p,\infty} (\Psi_1 \otimes \Psi_2)$ or, alternatively,  if and only if $k\in \mathfrak{l}^{\, p,\infty}_{1/(w_1\otimes w_2)}(I\times J)$; and in that case
\begin{equation}
\| O \|_{\H_{w_1}^1 (\Psi_1) \to \H_{1/w_2}^p (\Psi_2)}
\asymp \| K \|_{\mathfrak{H}_{1/(w_1 \otimes w_2)}^{p,\infty} (\Psi_1 \otimes \Psi_2)}=\|k\|_{\mathfrak{l}^{\, p,\infty}_{1/(w_1\otimes w_2)}(I\times J)},
\label{eq:Supergeiler_satz_Normequivalenz_1}
\end{equation}
 \item\label{thm:super-ii} A linear operator $O:\H_{w_1}^{p} (\Psi_1)\to\H_{1/w_2}^{\infty} (\Psi_2)$ is bounded  if and only if its kernel $K$ is contained in  $\H_{1/(w_1 \otimes w_2)}^{q, \infty} (\Psi_1 \otimes \Psi_2)$ or, alternatively, if and only if $k\in \ell^{ p,\infty}_{1/(w_1 \otimes w_2)}(I\times J)$; and in that case
\begin{equation}
\| O \|_{\H_{w_1}^{p} (\Psi_1) \to \H_{1/w_2}^{\infty} (\Psi_2)}
\asymp \| K \|_{\H_{1/(w_1 \otimes w_2)}^{q,  \infty} (\Psi_1 \otimes \Psi_2)} = \|k\|_{\ell^{ p,\infty}_{1/(w_1 \otimes w_2)}(I\times J)}.
\label{eq:Supergeiler_satz_Normequivalenz_2}
\end{equation}
\end{enumerate}
\end{theorem}
%\textcolor{blue}{**** New shorter proof***}
\noindent \proof Ad   \ref{thm:super-i}: %From \eqref{eq:Inklusion} we infer that $H_{w_1}^1 (\Psi_1) \subseteq H_{m_1}^1 (\Psi_1)$ as $m_1$ is supposed to be $w_1$-moderate and that $H_{m_2}^{\, p_1} (\Psi_2) \subseteq H_{1/w_2}^{\infty} (\Psi_2)$ as $m_2$ is supposed to be $w_2$-moderate. Therefore we are in a position to apply the kernel to the operator $O$; 
The fact that $O:\H_{w_1}^1 (\Psi_1)\to\H_{1/w_2}^p (\Psi_2)$ is bounded if and only if  $k\in \mathfrak{l}^{\, p,\infty}_{1/(w_1\otimes w_2)}(I\times J)$ was established in \cite[Proposition~6,~eq.~(24)]{Balazs_2017}. To show the second equivalence, we first observe that if $O:\H_{w_1}^1 (\Psi_1)\to\H_{1/w_2}^p (\Psi_2)$ is bounded, then
$$
\|Of\|_{\H^\infty_{1/w_2}(\Psi_2)}\le 
\|Of\|_{\H^p_{1/w_2}(\Psi_2)}\lesssim \|f\|_{\H^1_{w_1}(\Psi_1)}.
$$
According to Theorem~\ref{thm3}, there exists a unique kernel $K \in \H_{1 / (w_1 \otimes w_2)}^{\infty,  \infty} \left( \Psi_1 \otimes \Psi_2 \right)$ that satisfies \eqref{eq:Correspondence_1}. In particular,
$$
\langle K,\widetilde{\psi}_{1,\, i}\otimes \widetilde{\psi}_{2,\, j}\rangle_{H^{\infty,\infty}_{1/(w_1\otimes w_2)}(\Psi_1\otimes \Psi_2),H^{1,1}_{w_1\otimes w_2}(\Psi_1\otimes \Psi_2)}\hspace{-1pt}=\hspace{-1pt}\big\langle O\widetilde{\psi}_{1,\, i},\widetilde{\psi}_{2,\, j}\big\rangle_{\H^\infty_{1/w_2}(\Psi_2),\H_{w_2}^1(\Psi_2)}\hspace{-1pt}=\hspace{-1pt}k_{i,j},
$$
which shows that $\| K \|_{\mathfrak{H}_{1/(w_1 \otimes w_2)}^{p,\infty} (\Psi_1 \otimes \Psi_2)}=\|k\|_{\mathfrak{l}^{\, p,\infty}_{1/(w_1\otimes w_2)}(I\times J)}<\infty$ if $O$ is bounded. 

If, on the other hand, $K\in \mathfrak{H}_{1/(w_1 \otimes w_2)}^{p,\infty} (\Psi_1 \otimes \Psi_2)\subset \H_{1/(w_1 \otimes w_2)}^{\infty,\infty} (\Psi_1 \otimes \Psi_2)$, then from Theorem~\ref{kernsatz} we infer that there exists a bounded operator $O:\H^1_{w_1}(\Psi_1)\to\H^\infty_{1/w_2}(\Psi_2)$  satisfying \eqref{eq:Correspondence_1}. Therefore, we conclude that $\| K \|_{\mathfrak{H}_{1/(w_1 \otimes w_2)}^{p,\infty} (\Psi_1 \otimes \Psi_2)}=\|k\|_{\mathfrak{l}^{\, p,\infty}_{1/(w_1\otimes w_2)}(I\times J)}$
which  implies the boundedness of $O:\H^1_{w_1}(\Psi_1)\to\H^p_{1/w_2}(\Psi_2)$.

Ad \ref{thm:super-ii}: The equivalence of the boundedness of the operator and the property of its Galerkin matrix is not explicitly stated in \cite[Proposition~6]{Balazs_2017}. This can however be done by essentially using the same arguments as those used in \cite[Proposition~6]{Balazs_2017} and applying \cite[Proposition~2.4]{tao}.
The rest of the proof is similar to the one of  \ref{thm:super-i} with the only difference that here we need to make use of the inequalities
$$
\|Of\|_{\H^\infty_{1/w_2}(\Psi_2)}\lesssim  
\|f\|_{\H^p_{w_1}(\Psi_1)}\le  \|f\|_{\H^1_{w_1}(\Psi_1)},
$$
which ensures the existence of a kernel $K$.\pbox

This result allows us to prove that the co-orbit spaces associated with the mixed norm spaces of Theorem~\ref{Supergeiler_Satz} do not depend on the particular tensor product of  \lf frames as long as all involved cross-Gram matrices belong to $\mathcal{A}_n$. 

\begin{corollary} \label{cor}
 Let $1\le p \le \infty$, $\Psi_n,\Phi_n\subset \H_n$ be   $\mathcal{A}_n$-localised frames, and $w_n$ an $\mathcal{A}_n$-admissible weight $n=1,  2$. If $G_{\Phi_n,\widetilde{\Psi}_n}\in\mathcal{A}_n$ and $G_{\Psi_n,\widetilde{\Phi}_n}\in\mathcal{A}_n$, $n=1,2$, then   $\H_{w_1 \otimes w_2}^{p,\infty} \left( \Psi_1 \otimes \Psi_2 \right)=\H_{w_1 \otimes w_2}^{p,\infty} \left( \Phi_1 \otimes \Phi_2 \right)$ and  $\mathfrak{H}_{w_1 \otimes w_2}^{p,\infty} \left( \Psi_1 \otimes \Psi_2 \right)=\mathfrak{H}_{w_1 \otimes w_2}^{p,\infty} \left( \Phi_1 \otimes \Phi_2 \right)$, with their norms being equivalent.
\end{corollary}
\proof
We shall only prove the former identity as the latter follows from a similar argument. Let $K_\Psi\in \H_{w_1 \otimes w_2}^{p,\infty} \left( \Psi_1 \otimes \Psi_2 \right)$ and $O:\H_{1/w_1}^{p} (\Psi_1)\to\H_{1/w_2}^{\infty} (\Psi_2)$ be the corresponding bounded operator according to Theorem~\ref{Supergeiler_Satz}~\ref{thm:super-ii}.  From Proposition~\ref{prop:properties}~\ref{en:ii} we know that $\H^{p}_{1/ w_i}(\Psi_i)=\H^{p}_{1/ w_i}(\Phi_i)$ and thus $O:\H_{1/w_1}^{p} (\Phi_1)\to\H_{1/w_2}^{\infty} (\Phi_2)$ is bounded with equivalent operator norm. Applying Theorem~\ref{Supergeiler_Satz}~\ref{thm:super-ii} again, we find an element $K_\Phi\in \H_{w_1 \otimes w_2}^{p,\infty} \left( \Phi_1 \otimes \Phi_2 \right)$   such that 
\begin{align*}
\|K_\Psi\|_{\H_{w_1 \otimes w_2}^{p,\infty} \left( \Psi_1 \otimes \Psi_2 \right)}&\asymp 
\|O\|_{\H_{w_1}^{p} (\Psi_1)\to\H_{w_2}^{\infty} (\Psi_2) }
\\ &\asymp 
\|O\|_{\H_{w_1}^{p} (\Phi_1)\to\H_{w_2}^{\infty} (\Phi_2) }\asymp\|K_\Phi\|_{\H_{w_1 \otimes w_2}^{p,\infty} \left( \Phi_1 \otimes \Phi_2 \right)}. 
\end{align*}
It remains to show that $K_\Psi=K_\Phi$. This, however, follows from the uniqueness of the kernel in Theorem~\ref{kernsatz} %($O$ is also bounded as an operator mapping $\H_{1/w_1}^{1} (\Phi_1)$ on $\H_{w_2}^{\infty} (\Phi_2)$
and therefore $\H^{\infty,\infty}_{w_1\otimes w_2 }(\Psi_1\otimes\Psi_2) = \H^{\infty,\infty}_{w_1\otimes w_2 }(\Phi_1\otimes\Phi_2)$ as long as the cross-Gram matrices of $\Psi_1$ and $\Phi_1$, on the one hand, and $\Psi_2$ and $\Phi_2$, on the other hand, are contained in $\mathcal{A}_1$ and $\mathcal{A}_2$ respectively.
\pbox

Finally, we give a sufficient condition on the kernel to ensure that the corresponding operator is contained in the Schatten-$p$ class.
\begin{proposition}
    Let $\mathcal{S}_p(\H_1,\H_2)$ be  the Schatten-$p$ class of operators mapping $\H_1$ to $\H_2$, $1\le p\le 2$, $\Psi_n$ be $\mathcal{A}_n$-localised frames for $\H_n$ for $n=1,  2$, and $O:\H_1\to\H_2$ a bounded operator with a kernel $K$.
    If $K\in\mathfrak{H}^{2,p}(\Psi_1\otimes\Psi_2)$, then $O\in \mathcal{S}_p(\H_1,\H_2)$. 
    In particular, if     $K\in\H^{p,p}(\Psi_1\otimes\Psi_2)$,  then $O\in \mathcal{S}_p(\H_1,\H_2)$.
\end{proposition}
\proof From
\cite[Theorem~17]{schatten} we know that, if $\big\{\|Of_i\|_{\H_2}\big\}_{i\in I}\in \ell^p(I)$ for some frame $\{f_i\}_{i\in I}\subset\H_1$, then $O\in \mathcal{S}_p(\H_1,\H_2)$. Using $\Psi_1$ as such a frame results in
\begin{align*}
\sum_{i\in I}\big\|O\widetilde{\psi}_{1,\, i}\big\|_{\H_2}^p&\asymp \sum_{i\in I}\Big(\sum_{j\in J}\big|\big\langle O\widetilde{\psi}_{1,\, i},\widetilde{\psi}_{2,\, j}\big\rangle_{\H_2}\big|^2\Big)^{p/2}
\\
&= \sum_{i\in I}\Big(\sum_{j\in J}\big|\big\langle K,\widetilde{\psi}_{1,\, i}\otimes \widetilde{\psi}_{2,\, j}\big\rangle_{\H_1\otimes \H_2}\big|^2\Big)^{p/2}=\|K\|_{\mathfrak{H}^{2,1}(\Psi_1\otimes\Psi_2)}.
\end{align*}
The second statement follows immediately once we   observe that $\|K\|_{\mathfrak{H}^{2,1}(\Psi_1\otimes\Psi_2)}\le\|K\|_{\H^{p,p}(\Psi_1\otimes\Psi_2)} $.
\pbox

\section*{Acknowledgements}

\noindent We would like to thank Hans Georg Feichtinger for fruitful discussions on the contents of this article. Dimitri Bytchenkoff is grateful to Patrik Wahlberg of Polytechnic of Turin for a helpful conversation on kernel theorems.
This research was funded by the Austrian Science Fund (FWF) 10.55776/P34624 (P.B. and D.B.) and 10.55776/Y1199 (M.S.).
\begin{comment}
For open
access purposes, the authors have applied a CC BY public copyright license to any author-accepted manuscript
version arising from this submission.
\end{comment}


\begin{thebibliography}{99}

\bibitem{Balazs_2008}
  P. Balazs,
  \textit{Matrix representation of operators using frames}.
  Sampl. Theory Signal Image Process. 7, (2008) 39-54.

\bibitem{Balazs_2008_bis}
  P. Balazs,
  \textit{Hilbert-Schmidt operators and frames - classification, best approximation by multipliers and algorithms}.
  Int. J. Wavelets Multiresolution Inf. Process. 6(2), (2008) 315-330.

%\bibitem{Balazs_2011} \textcolor{red}{  P. Balazs, D. Stoeva and J.-P. Antoine, \textit{Classification of general sequences by frame-related operators}.  Sampl. Theory Signal Image Process. 10, (2011) 151-170.}

\bibitem{Balazs_2017}
  P. Balazs and K. Gröchenig,
  \textit{A Guide to localized frames and applications to Galerkin-like representations of operators, in I. Pesenson, Q. Le Gia, A. Mayeli, H. Mhaskar, D.-X. Zhou, Frames and other bases in abstract and function spaces}.
  Birkhäuser, Cham (2017).

\bibitem{Balazs_2019}
  P. Balazs, K. Gröchenig and M. Speckbacher,
  \textit{Kernel theorems in coorbit theory}.
  Trans. Amer. Math. Soc. Ser. B 6, (2019) 346-364.

 \bibitem{Besov_1961}
  O. V. Besov
  \textit{On a family of function spaces in connection with embeddings and extensions}. Tr. Mat. Inst. Steklova 60 (1961), 42–81.


  \bibitem{Borup_2007}
 L. Borup and M. Nielsen,
  \textit{Frame decomposition of decomposition spaces}.
  J. Fourier Anal. Appl. 13(1), (2007) 39-70.

   \bibitem{Bytchenkoff_2020}
 D. Bytchenkoff and F. Voigtlaender,
  \textit{Design and properties of wave packet smoothness spaces}.
  J. Math. Pures Appl. 133, (2020) 185-262.

  \bibitem{Bytchenkoff_2021}
 D. Bytchenkoff,
  \textit{Construction of Banach frames and atomic decompositions
of anisotropic Besov spaces}.
  Ann. Mat. Pura Appl. 200, (2021) 1341-1365.

  \bibitem{sampta} D. Bytchenkoff, M. Speckbacher and P. Balazs, \textit{Outer kernel theorem for co-orbit spaces of localised frames}. In: Proceedings of SampTA 2023, (2023). 

  \bibitem{schatten} H.  Bingyang, L. H. Khoi and K. Zhu, \textit{ Frames and operators in {S}chatten classes}. Houston J.  Math. 41(4), (2013) 1191-1219.

\bibitem{Casazza_2013}
   P. G. Casazza and G. Kutyniok
  \textit{Finite frames theory and applications}. Applied and numerical
harmonic analysis. Boston, MA
  Birkhäuser, Boston (2013).    

\bibitem{Christensen_2008}
  O. Christensen,
  \textit{Introduction to Frames and Riesz Bases}.
  Birkhäuser, Boston (2008).  

 \bibitem{conway} J. B. Conway, \textit{A Course in Functional Analysis}.  Graduate Texts in Mathematics,
vol. 96, Springer-Verlag, New York, (1990).

\bibitem{Cordero_2019}
  E. Cordero and F. Nicola,
  \textit{Kernel theorems for modulation spaces}.
  J. Fourier Anal. Appl. 25, (2019) 131-144.
  
%\bibitem{Daubechies_1992} I. Daubechies, \textit{Ten lectures on wavelets}, Society For Industrial and Applied Mathematics (1992).

%\bibitem{Dieudonne_1950}  J. Dieudonné,  \textit{Natural homomorphisms in Banach spaces}.  Proc. Amer. Math. Soc. MR 11, 524, (1950) 54-59.   
 
\bibitem{Duffin_1952}
  R. J. Duffin and A. C. Schaeffer,
  \textit{A class of nonharmonic Fourier series}.
  Trans. Amer. Math. Soc. 72, (1952) 341-154.    

\bibitem{Feichtinger_1980}
  H. G. Feichtinger,
  \textit{Un espace de Banach de distributions tempérées sur les groupes localement compacts abéliens}.
  C. R. Acad. Sci. Paris Sér. A-B 290(17), (1980) A791-A794.

\bibitem{Feichtinger_1981}
  H. G. Feichtinger,
  \textit{On a new Segal algebra}.
  Monatsh. Math. 92, (1981) 269-289.
  
\bibitem{Feichtinger_1985}  H. G. Feichtinger and P. Gröbner, \textit{Banach spaces of distributions defined by decomposition method}. Math. Nachr. 123, (1985) 97-120.

\bibitem{Feichtinger_1988}
  H. G. Feichtinger and K. Gröchenig,
  \textit{A unified approach to atomic decomposition via integrable group representations, In : Proc. Conference on Functions, Spaces and Applications, Lund 1986},
  Springer Lect. Notes Math. (1988) 52-73.

\bibitem{Feichtinger_1989} H. G. Feichtinger and K. Gröchenig, \textit{Banach spaces related to integrable group representations and their atomic decompositions, I.}.  J. Funct. Anal. 86(2), (1989) 307-340.  

\bibitem{Feichtinger_1989_2}  H. G. Feichtinger and K. Gröchenig,  \textit{Banach spaces related to integrable group representations and their atomic decompositions. II.}.  Monatsh. Math. 108(2-3), (1989) 129-148. 

  \bibitem{Feichtinger_Jakobsen}
  H. G. Feichtinger and M. S. Jakobsen,
  \textit{The inner kernel theorem for a certain {S}egal algebra}.
  Monatsh. Math. 198, (2022) 675–715. 

 % \bibitem{Feichtinger-Kozek} H. G. Feichtinger and W. Kozek, \textit{Quantization of TF lattice-invariant operators on elementary LCA groups}. In: H. G. Feichtinger, T. Strohmer  (eds.) Gabor Analysis and Algorithms. Applied and Numerical Harmonic Analysis. Birkhäuser, Boston (1998).

\bibitem{Fornasier_2005}
  M. Fornasier and K. Gröchenig, 
  \textit{Intrinsic localization of frames}.
   Constr. Approx. 22, (2005) 395-415.   

   \bibitem{Fornasier_Rauhut}
  M. Fornasier and H. Rauhut, 
  \textit{Continuous frames, function spaces, and the discretization problem}.
  J. Fourier Anal. Appl. 11(3), (2005) 245-287.  

   \bibitem{Futamura_2009}
  F. Futamura, 
  \textit{Localizable operators and the construction of localized frames}.
  Proc. Amer. Math. Soc. 137(12), (2009) 4187-4197.    

\bibitem{Galerkin_1915}
  B. G. Galerkin,
  \textit{Rods and plates: series occurring in various questions concerning the elastic equilibrium of rods and plates}.
  Vestnik Inzhenerov i Tekhnikov 19, (1915) 897-908. 

  \bibitem{geshi} I. M. Gelfand and G. E. Shilov, \textit{Generalized Functions. Vol. 2, Spaces of fundamental and
generalized functions}. Translated from the 1958 Russian original,  AMS Chelsea Publishing, Providence, RI, (2016).
  
%\bibitem{Groebner_1992}  P. Gröbner,  \textit{Banachräume glatter Funktionen und Zerlegungsmethoden}.  DPhil thesis, University of Vienna (1992).           

\bibitem{Groechenig_1991}   K. Gröchenig, \textit{Describing functions: atomic decompositions versus frames}. Monatsh. Math. 112, (1991) 1-41.
  
%\bibitem{Groechenig_2001}   K. Gröchenig, \textit{Foundations of time-frequency analysis},  Birkhäuser, Basel (2001).
  
\bibitem{Groechenig_2003}
  K. Gröchenig,
  \textit{Localized frames are finite unions of Riesz sequences}.
   Adv. Comp. Math. 18, (2003)  149-157.

\bibitem{Groechenig_2004}
  K. Gröchenig,
  \textit{Localization of frames, Banach frames and the invertibility of the frame operator}.
   J. Fourier Anal. Appl. 10(2) (2004)  105-132.   

%\bibitem{Groechenig_2006}     K. Gröchenig,  \textit{A pedestrian’s approach to pseudodifferential operators, In : Harmonic analysis and applications}.   Appl. Numer. Harmon. Anal., Birkhäuser, Boston (2006), 139-169. 

   \bibitem{Schur} K. Gr\"ochenig and M. Leinert, \textit{Symmetry of matrix algebras and symbolic calculus for infinite
matrices}. Trans. Amer. Math. Soc.  358, (2006), 2695-2711


\bibitem{testfunction} L. H\"ormander, \textit{The Analysis of Linear partial differential operators. I,  Distribution theory
and Fourier analysis}.   Springer, Berlin (1990).
  
\bibitem{Jaffard_1990}
  S. Jaffard,
  \textit{Propriétés des matrices "bien localisées" près de leur diagonale et quelquels applications}.
   Ann. Inst. Henri Poincar{\'e} 7, (1990), 461-476.

   \bibitem{Losert} V. Losert, \textit{A characterization of the minimal strongly character invariant {S}egal algebra}. Ann. Inst. Fourier 30, (1980) 129-139.

   \bibitem{maddox} I. J. Maddox, \textit{Infinite Matrices of Operators}. Lecture Notes in Mathematics, vol. 786, Springer, Berlin (1980).

  \bibitem{Nielsen_2012}
  M. Nielsen and K. N. Rasmussen,
  \textit{Compactly supported frames for decomposition spaces}.
   J. Fourier Anal. Appl. 18(1), (2012), 87-117. 

\bibitem{Nielsen_2014}
  M. Nielsen,
  \textit{Frames for decomposition spaces generated by a single function}.
   Collect. Math. 65, (2014), 183-201. 
      \bibitem{Ryan_2002} R. A. Ryan, \textit{Introduction to Tensor Products of Banach Spaces}, Springer, London (2002).

   %\bibitem{megg} R. Megginson, \textit{An Introduction to Banach Space Theory}. Springer, New York (1998).

   %\bibitem{Schatten_1960} R. Schatten, \textit{Norm ideals of completely continuous operators}. Springer, Berlin (1960).

   \bibitem{Sjöstrand} J. Sj\"ostrand, \textit{Wiener-type algebras of pseudodifferential operators}. Centre de Mathematiques,
Ecole Polytechnique. Palaiseau France, S\'eminaire 1994–1995, 
(1994).

\bibitem{tao} T. Tao, \textit{Lecture Notes 2 for 247 A: Fourier Analysis}. \url{www.math.ucla.edu/~tao/247a.1.06f/notes2.pdf}.

\bibitem{Voigtlaender _2022}
  F. Voigtlaender,
  \textit{Structured, compactly supported Banach frame decompositions of decomposition spaces}.
   Diss. Math. 575, (2022), 1-179. 

   %\bibitem{Werner_1995} D. Werner, \textit{Funktionalanalysis}. Springer, Berlin (1995).
   
\end{thebibliography}
\end{document}